\documentclass[reqno,12pt]{amsart}

\usepackage{amsmath}

\usepackage{amsfonts}

\usepackage{amssymb}

\usepackage{eufrak}

\usepackage{amscd}
\usepackage{amsthm}
\usepackage{amstext}
\usepackage[all]{xy}
   \newcommand{\Hom}{\operatorname{Hom}}
   
 \newcommand{\Ext}{\operatorname{Ext}}
\newcommand{\Ad}{\operatorname{Ad}}

\newcommand{\id}{\operatorname{id}}

\newcommand{\im}{\operatorname{im}}
\newcommand{\diag}{\operatorname{diag}}

\newcommand{\ess}{\operatorname{ess}}
\newcommand{\re}{\operatorname{Re}}

    \newcommand{\ev}{{ev}}

   \theoremstyle{plain}%default
   \newtheorem{thm}{Theorem}[section]
   \newtheorem{prop}[thm]{Proposition}
   \newtheorem{lemma}[thm]{Lemma}
   
   \theoremstyle{definition}

   \theoremstyle{remark}
   
   \newtheorem{remark}[thm]{Remark}

   \numberwithin{equation}{section}

\title[Relative K-homology]{Relative K-homology and normal operators}
\author{Vladimir Manuilov and Klaus Thomsen}

        \date{\today}
\address{Department of Mechanics and Mathematics, Moscow State
University, Leninskie Gory, Moscow, 119992, Russia}
\email{manuilov@mech.math.msu.su}

  \address{IMF, Department of Mathematics, Ny Munkegade, 8000 Aarhus C, Denmark}

\email{matkt@mi.aau.dk}

\begin{document}

        %\thanks{My warmest thanks to my mom}
        %\email{matkt@mi.aau.dk}
        %\address{Department of Mathematics, Ny Munkegade, 8000 Aarhus C, Denmark}
        %\dedicatory{Dedicated to S\o ren Have Hansen on the occasion of his birth}
        %\keywords{TeX, Mathematics}
        %\subjclass{14C10, 93D20}
        %\commby{Jens Peter S\ae r}

%Slut paa  'amsart' style specifikke ting.

%Flg. kommando kan justere hvor detaljeret Indholdsfortegnelse (Table Of Contents)
%bliver.
%\setcounter{tocdepth}{3}

\maketitle

%\tableofcontents

%----------Selve teksten -------------------------------------------------------

\section{Introduction.}

Let $X$ be a compact metric space. By results of Brown, Douglas
and Fillmore, \cite{BDF2}, the $K$-homology of $X$ is realized by
$\Ext(X)$, the equivalence classes of unital and essential
extensions of $C(X)$ by the compact operators $\mathbb K$ on a
separable infinite dimensional Hilbert space $H$, or equivalently,
the equivalence classes of unital and injective $*$-homomorphisms
$C(X) \to Q$, where $Q = \mathbb L(H)/\mathbb K$ is the Calkin
algebra. This discovery came out of questions and problems related
to essential normal operators, and it led quickly to the
development of a vast new area of mathematics which combines
operator theory with algebraic topology. In particular, the
BDF-theory was generalized by Kasparov in form of KK-theory, which
has proven to be a powerful tool in the theory of operator
algebras as well as in algebraic topology.

It is the purpose of the present paper to develop a relative
theory in this context. The point of departure here is a
generalization of the six-term exact sequence of extension theory
which relates the group of extensions of a unital $C^*$-algebra to
the group of unital extensions. This sequence was discovered by
Skandalis, cf. \cite{S}, and a construction of it was presented in
\cite{MT}. It is the latter construction which we here generalize
to get a relative extension theory. Subsequently we investigate
the relative K-homology which arises from it by specializing to
abelian $C^*$-algebras. It turns out that relative K-homology
carries substantial information also in the operator theoretic
setting from which the BDF-theory was developed, cf. \cite{BDF1}, and
we conclude the paper by extracting some of this information.

In the remaining part of this introduction we give a more detailed
account of the content of the paper. Let $A$ be a $C^*$-algebra,
$J \subseteq A$ a $C^*$-subalgebra, and let $B$ be a stable
$C^*$-algebra. Under modest assumptions we organize the
$C^*$-extensions of $A$ by $B$ that are trivial when restricted
onto $J$ to become a semi-group $\Ext_J(A,B)$ which is the
semi-group $\Ext(A,B)$ of Kasparov, \cite{K1}, when $J = \{0\}$.
The group $\Ext_J^{-1}(A,B)$ of invertible elements in
$\Ext_J(A,B)$ can be effectively computed by a six-term exact
sequence which generalizes the excision six-term exact sequence in
the first variable of $KK$-theory, and it turns out that there is
a natural identification $\Ext^{-1}_J(A,B) =
KK\left(C_i,B\right)$, where $C_i$ is the mapping cone of the
inclusion $i : J \to A$. Thus, as an abstract group, the relative
extension group is a familiar object, and the six-term exact
sequence which calculates it, is a version of the Puppe exact
sequence of Cuntz and Skandalis; \cite{CS}. But the realization of
$KK(C_i,B)$ as a relative extension group has non-trivial
consequences already in the set-up from which KK-theory developed,
namely the setting of (essential) normal operators, and the second
half of the paper is devoted to the extraction of the information
which the relative extension group contains about normal operators
when specialized to the case where $B = \mathbb K$ and $X$ and $Y$
are compact metric spaces, and $f : X \to Y$ is a continuous
surjection giving rise to an embedding of $J = C(Y)$ into $A
=C(X)$. In this setting $\Ext_J(A,\mathbb K)$ is a group, and we
denote it by $\Ext_{Y,f}(X)$. As an abstract group this is the
even K-homology of the mapping cone of $f$, and the above
mentioned six-term exact sequence takes the form
\begin{equation}\label{introsix}
\begin{xymatrix}{
 \Ext_{Y,f}(X)  \ar[r]  &   K_1(X)  \ar[r]^-{f_*} & K_1(Y)  \ar[d] \\
K_0(Y)  \ar[u] &   K_0(X) \ar[l]^-{f_*} &  \Ext_{SY,Sf}(SX) \ar[l] }
\end{xymatrix}
\end{equation}
where $S$ is the reduced suspension. An element of $\Ext_{Y,f}(X)$
consists of a commuting diagram
\begin{equation}\label{intro}
\begin{xymatrix}{
C(X) \ar[r]^-{\varphi}  &  Q \\
C(Y) \ar[u]^-{f^*}  \ar[r]^-{\varphi_0}  & {\mathbb L(H)} \ar[u]  }
\end{xymatrix}
\end{equation}
where $\varphi$ and $\varphi_0$ are unital and injective
$*$-homomorphisms. Thus $\varphi$ is an extension of $C(X)$ by
$\mathbb K$, in the sense of Brown, Douglas and Fillmore, which is
trivial (or split) when restricted to $C(Y)$, and $\varphi_0$ is a
specified splitting of the restriction. $\Ext_{Y,f}(X)$ can be
defined as the homotopy classes of such diagrams, or pairs
$(\varphi,\varphi_0)$, but as one would expect from experience
with BDF-theory and KK-theory, the group admits several other
descriptions where the equivalence relation is seemingly stronger
and/or the diagrams are required to have special properties. In
particular, triviality of the diagram (\ref{intro}) is equivalent
to the existence of $*$-homomorphisms $\psi_n : C(X) \to \mathbb
L(H)$ such that the upper triangle in the diagram
\begin{equation}\label{intro2}
\begin{xymatrix}{
C(X) \ar[r]^-{\varphi} \ar[dr]^-{\psi_n} &  Q \\
C(Y) \ar[u]^-{f^*}  \ar[r]^-{\varphi_0}  & {\mathbb L(H)} \ar[u]  }
\end{xymatrix}
\end{equation}
commutes for each $n$, and the lower triangle asymptotically
commutes in the sense that $\lim_{n \to \infty} \psi_n \circ
f^*(g) = \varphi_0(g)$ for all $g \in C(Y)$. Thus the relative
extension group $\Ext_{Y,f}(X)$ presents the obstructions for the
existence of a splitting of the whole extension $\varphi$ which
respects the given splitting over $C(Y)$ up to any given
tolerance. These obstructions are naturally divided in two
classes, where the first is the rather obvious obstruction that
the diagram (\ref{intro2}) can only exist when the extension
$\varphi$ is split. This obstruction is described by the presence
of an obvious map $\Ext_{Y,f}(X) \to  \Ext(X)$. In many cases this
map is injective, and then the obvious obstruction is the only
obstruction. But generally the map to $\Ext(X)$ is not injective,
and the kernel of it consists of the non-trivial obstructions --
those that arise because we insist that the given splitting over
$C(Y)$ should be respected, at least asymptotically. The six-term
exact sequence (\ref{introsix}) shows that the kernel of the
forgetful map $\Ext_{Y,f}(X) \to  \Ext(X)$ is isomorphic to the
co-kernel of the map $f_* : K_0(X) \to K_0(Y)$. This part of the
relative K-homology contains the obstructions for finding a
$*$-homomorphic lift $C(X) \to \mathbb L(H)$ of $\varphi$ which
agrees with $\varphi_0$ on $C(Y)$ up to an arbitrarily small
compact perturbation. We show that this part of the relative
K-homology vanishes in many cases, and in particular when $Y$ is a
compact subset of the complex plane $\mathbb C$. This then serves
as the main ingredience in the proof of the following
operator-theoretic fact:

\begin{thm}\label{introthm} Let $M,N_1,N_2, N_3, \dots, N_k$ be
bounded  normal operators such that $N_iN_j = N_jN_i$ for all
$i,j$, and let $F$ be a continuous function from the joint
spectrum of the $N_i$'s onto the spectrum of $M$ such that
 $$
F\left(N_1,N_2, \dots, N_k\right) - M \in \mathbb K .
 $$
For every $\epsilon > 0$ there are normal operators
$N_1^{\epsilon}, N_2^{\epsilon}, N_3^{\epsilon}, \dots,
N_k^{\epsilon} $ such that $N_i^{\epsilon}{N_j^{\epsilon}} =
{N_j^{\epsilon}}N_i^{\epsilon}$, $N_i-N_i^{\epsilon} \in \mathbb
K$ for all $i,j$, and
 $$
\left\|F\left(N_1^{\epsilon},N_2^{\epsilon}, \dots, N_k^{\epsilon}\right) - M \right\| \leq \epsilon.
 $$
\end{thm}

Examples show that the $\epsilon > 0$ can not in general be replaced by $0$.

\section{The relative extension group.} \label{sec2}

We begin by recalling the definition of the group of
$C^*$-extensions, as it appears in KK-theory. Let $A,B$ be
separable $C^*$-algebras, $B$ stable. As is well-known, the
$C^*$-algebra extensions of $A$ by $B$ can be identified with
$\Hom (A,Q(B))$, the set of $*$-homomorphisms $A \to Q(B)$, where
$Q(B) = M(B)/B$ is the generalized Calkin algebra. Let $q_B : M(B)
\to Q(B)$ be the quotient map. Two extensions $\varphi, \psi : A
\to Q(B)$ are \emph{unitarily equivalent} when there is a unitary
$u \in M(B)$ such that $\Ad q_B(u) \circ \psi = \varphi$. The
unitary equivalence classes of extensions of $A$ by $B$ have the
structure of an abelian semi-group thanks to the stability of $B$:
Choose isometries $V_1,V_2 \in M(B)$ such that $V_1V_1^* +
V_2V_2^* = 1$, and define the sum $\varphi \oplus \psi : A \to
Q(B)$ of $\varphi, \psi \in\Hom (A,Q(B))$ to be
\begin{equation}\label{B1}
(\psi \oplus \varphi)(a) = \Ad q_B(V_1) \circ \psi(a)     + \Ad q_B(V_2) \circ \varphi(a) .
\end{equation}
An extension $\varphi : A \to Q(B)$ is \emph{split} when there is
a $*$-homomorphism $\pi : A \to M(B)$ such that $\varphi = q_B
\circ \pi$. To trivialize the split extensions and obtain a
neutral element for the composition we declare two extensions
$\varphi, \psi : A \to Q(B)$ to be \emph{stably equivalent} when
there is a split extension $\pi$ such that $\psi \oplus \pi$ and
$\varphi \oplus \pi$ are unitarily equivalent. The semigroup of
stable equivalence classes of extensions of $A$ by $B$ is denoted
by $\Ext(A,B)$. As is well-documented by now, the semi-group is
generally not a group, and we denote by
 $$
\Ext^{-1}(A,B)
 $$
the abelian group of invertible elements in $\Ext(A,B)$.

An \emph{absorbing $*$-homomorphism} $\pi : A \to M(B)$ is a
$*$-homomorphism with the property that for  every completely
positive contraction $\varphi : A \to M(B)$ there is a sequence
$V_n \in M(B)$ of isometries such that $V_n^*\pi(a)V_n -
\varphi(a) \in B$ for all $n$, and $\lim_{n \to \infty}
V_n^*\pi(a)V_n = \varphi(a)$ for all $a \in A$. When $A$ is
unital, a \emph{unitally absorbing $*$-homomorphism} $\pi : A \to
M(B)$ is a unital $*$-homomorphism with the property that for
every completely positive contraction $\varphi : A \to B$ there is
a sequence $W_n \in M(B)$ such that $\lim_{n \to \infty}
W_n^*\pi(a)W_n = \varphi(a)$ for all $a \in A$, and $\lim_{n \to
\infty} W_n^*b = 0$ for all $b \in B$. We refer the reader to
\cite{Th1} for alternative characterizations of absorbing and
unitally absorbing $*$-homomorphisms which justify the names, and
a proof that they always exist in the separable case. Of
particular importance here is the \emph{essential uniqueness} of
such $*$-homomorphisms. Specifically, when $\pi,\lambda : A \to
M(B)$ are $*$-homomorphisms that are either both absorbing or both
unitally absorbing, there is a sequence $U_n$ of unitaries in
$M(B)$ such that $U_n\pi(a)U_n^* - \lambda(a) \in B$ for all $n$,
and $\lim_{n \to \infty} U_n\pi(a)U_n^* = \lambda(a)$ for all $a
\in A$.

Let now $J\subseteq A$ be a $C^*$-subalgebra of $A$, and consider
an absorbing $*$-homomorphism $\alpha_0 : A \to M(B)$. Set $\alpha
= q_B\circ \alpha_0 : A \to Q(B)$, and let
\begin{equation*}\label{B31}
\begin{xymatrix}{
0 \ar[r]  & B \ar[r] & E_0 \ar[r]  &  J \ar[r]  & 0}
\end{xymatrix}
\end{equation*}
be the extension of $J$ by $B$ whose Busby invariant is
$\alpha|_J$. Let $i : J \to A$ be the inclusion. We consider
extensions $E$ of $A$ by $B$ which fit into a commutative diagram
\begin{equation*}\label{B2}
\begin{xymatrix}{
0 \ar[r]  & B \ar@{=}[d] \ar[r] & E \ar[r]  &  A \ar[r]  & 0\\
0 \ar[r]  & B  \ar[r]   &  E_0 \ar[r]  \ar[u]  & J \ar[u]_-i \ar[r]  & 0 }
\end{xymatrix}
\end{equation*}
of $C^*$-algebras. In terms of the Busby invariant this
corresponds to extensions $\varphi : A \to Q(B)$ such that
$\varphi|_J = \alpha|_J$. We say that $\varphi$ \emph{equals
$\alpha$ on $J$}. Two such extensions, $\varphi, \psi : A \to
Q(B)$, that both equal $\alpha$ on $J$, are said to be
\emph{unitarily equivalent} when there is a unitary $v$ connected
to $1$ in the unitary group of the relative commutant $\alpha(J)'
\cap Q(B)$ such that $\Ad v \circ \varphi = \psi$.

\begin{lemma}\label{B91} For each $n \in \mathbb N$, there are
isometries $v_1,v_2, \dots, v_n$ in $\alpha(A)' \cap Q(B)$ such
that $v_i^*v_j = 0,i \neq j$, and $\sum_{i=1}^n v_iv_i^* = 1$.
\begin{proof}
Let $S_1,S_2,S_3, \dots$ be a sequence of isometries in $M(B)$
such that $S_i^*S_j = 0, i \neq j$, and $\sum_{i=1}^{\infty}
S_iS_i^* = 1$ with convergence in the strict topology. Set
$\beta_0(a) = \sum_{i=1}^{\infty} S_i\alpha_0(a)S_i^*$ for all $a
\in A$, and note that $\beta_0$ is absorbing because $\alpha_0$
is. Set $W_i = \sum_{j=1}^{\infty} S_{i + jn}S_j^*$. Then
$W_1,W_2, \dots, W_n$ are isometries in $ \beta_0(A)' \cap M(B)$
such that $W_i^*W_j = 0$ when $i \neq j$, and $\sum_{i=1}^n
W_iW_i^* = 1$. The essential uniqueness property of absorbing
$*$-homomorphisms guarantees the existence of a unitary $U \in
M(B)$ such that $\Ad U \circ \beta_0 (a) - \alpha_0(a) \in B$ for
all $a \in A$. Set $v_i = q_B\left( UW_iU^*\right), i = 1,2,
\dots, n$.
\end{proof}
\end{lemma}

Note that any choice of isometries $v_1, \dots, v_n$ as in Lemma
\ref{B91} gives us a $*$-isomorphism, $\Theta_n$, which maps the
$C^*$-algebra $M_n\left(\alpha(A)'\cap Q(B)\right)$ onto
$\alpha(A)'\cap Q(B)$. $\Theta_n$ is given by
\begin{equation*}\label{B93}
\Theta_n \left( \left(x_{ij}\right)\right) = \sum_{i,j =1}^n v_ix_{ij}v_j^* .
\end{equation*}
Any other choice of isometries as in Lemma \ref{B91} will result
in an isomorphism which is conjugate to $\Theta_n$ by a unitary
from $\alpha(A)' \cap Q(B)$. Thanks to Lemma \ref{B91} we can
define a composition $+_{\alpha}$ among the extensions of $A$ by
$B$ which agree with $\alpha$ on $J$:
\begin{equation}\label{B92}
\left(\varphi +_{\alpha} \psi\right)(a) =  w_1\varphi(a)w_1^* + w_2\psi(a)w_2^*,
\end{equation}
where $w_i \in \alpha(A)' \cap Q(B)$ are isometries such that
$w_1w_1^* + w_2w_2^* = 1$. To show that (\ref{B92}) gives the
unitary equivalence classes of extensions of $A$ by $B$ that agree
with $\alpha$ on $J$ the structure on an abelian semi-group,
consider the $*$-homomorphism $\beta_0$ introduced in the proof of
Lemma \ref{B91}, and set $\beta = q_B \circ \beta_0$. There are
then isometries $V_i, i = 1,2$, in $\beta_0(A)' \cap M(B)$ such
that $V_1V_1^* + V_2V_2^* = 1$. With these isometries we define
the sum $\varphi +_{\beta_0} \psi$ of two extensions, $\varphi,
\psi$, of $A$ by $B$ which both equal $\beta$ on $J$ to be
$\varphi +_{\beta_0} \psi = \Ad q_B(V_1) \circ \varphi + \Ad
q_B(V_2) \circ \psi $, in analogy with (\ref{B92}).

\begin{lemma}\label{B69}
Let $P \in \beta_0(A)'\cap M(B)$ be a projection such that both
$P$ and $1-P$ are Murray-von Neumann equivalent to $1$ in $
\beta_0(A)'\cap M(B)$. Let $U$ be a unitary in $\beta_0(A)'\cap
M(B)$ such that $UPU^* = P$. It follows that $U$ is connected to
$1$ in the unitary group of $ \beta_0(A)'\cap M(B)$.
\begin{proof}
Note first that $K_1\left(\beta_0(A)'\cap M(B)\right) = 0$ by
Lemma 3.1 of \cite{Th1}. Since $U = [PUP + (1-P)][P +
(1-P)U(1-P)]$, the lemma follows from this.
\end{proof}
\end{lemma}

\begin{lemma}\label{B73}
Let $\varphi, \psi, \lambda : A \to Q(B)$ be extensions of $A$ by
$B$ that both equal $\beta$ on $J$. There are then unitaries $S,T
\in \beta_0(A)'\cap M(B)$, connected to $1$ in the unitary group
of $ \beta_0(A)'\cap M(B),$ such that $\Ad q_B(S) \circ \left(
\varphi  +_{\beta_0} \psi\right) = \psi  +_{\beta_0} \varphi$ and
$\Ad q_B(T) \circ \left( \left( \varphi  +_{\beta_0} \psi\right)
+_{\beta_0} \lambda\right) = \varphi  +_{\beta_0} \left( \psi
+_{\beta_0} \lambda\right)$.
\begin{proof}
Let $\Theta: M_2(B) \to B$ be the $*$-isomorphism given by
 $$
\Theta \left( \begin{matrix} b_{11}  & b_{12} \\ b_{21}  & b_{22} \end{matrix} \right) = \sum_{i,j =1}^2 V_ib_{ij}V_j^* ,
 $$
and let $\overline{\Theta} : M(M_2(B)) = M_2(M(B) ) \to M(B)$ be
the $*$-isomorphism extending $\Theta$. Then $\overline{\Theta}
\circ \left( \begin{smallmatrix} \beta_0  & \\   & \beta_0
\end{smallmatrix} \right) = \beta_0$, so we see that
\begin{equation}\label{B70}
S = V_1V_2^* + V_2V_1^*  = \overline{\Theta}  \left( \begin{smallmatrix} 0  & 1 \\  1 & 0 \end{smallmatrix} \right)
\end{equation}
is connected to $1$ in the unitary group of $ \beta_0(A)'\cap
M(B)$. Since $\Ad q_B(S) \circ \left( \varphi  +_{\beta_0}
\psi\right) = \psi  +_{\beta_0} \varphi$, this proves the first
statement. To prove the second statement we identify $M_3(M(B))$
with $\mathbb L_B( B \oplus B \oplus B)$ -- the $C^*$-algebra of
adjointable operators on the Hilbert $B$-module $B^3$. Define
unitaries $W : B^3 \to B$ and $Z : B^3 \to B$ of Hilbert
$B$-modules such that $W(b_1,b_2,b_3) = V_1^2 b_1 + V_1V_2 b_2 +
V_2 b_3$ and $Z(b_1,b_2,b_3) = V_1b_1 + V_2V_1b_2 + V_2^2 b_3$.
Then $ZW^* \in M(B)$ and $\Ad q_B(ZW^*) \circ \left( \left(\varphi
+_{\beta_0} \psi\right)  +_{\beta_0} \chi\right) = \varphi
+_{\beta_0} \left( \psi  +_{\beta_0} \lambda\right)$. It remains
to show that $ZW^*$ is connected to $1$ in the unitary group of
$\beta_0(A)'\cap M(B)$. First observe that $ZW^*$ is connected to
 $$
Z  \left( \begin{smallmatrix}   &  & 1 \\   & 1  & \\ 1 &  &    \end{smallmatrix} \right)W^*
 $$
in the unitary group of $\beta_0(A)'\cap M(B)$. Since the unitary
$S$ from (\ref{B70}) is connected to $1$ in the unitary group of $
\beta_0(A)'\cap M(B)$ we see that $ZW^*$ is connected to
 $$
T = SZ  \left( \begin{smallmatrix}   &  & 1 \\   & 1  & \\ 1 &  &    \end{smallmatrix} \right)W^*
 $$
in the unitary group of $\beta_0(A)'\cap M(B)$. Note that
 $$
W^*V_2V_2^* W =  \left( \begin{smallmatrix} 0  &  &  \\   & 0  & \\  &  &  1  \end{smallmatrix} \right) \ \text{and}  \ Z^*V_1V_1^*Z = \left( \begin{smallmatrix} 1  &  &  \\   & 0  & \\  &  &  0  \end{smallmatrix} \right),
 $$
which implies that
 $$
Z  \left( \begin{smallmatrix}   &  & 1 \\   & 1  & \\ 1 &  &    \end{smallmatrix} \right)W^*V_2V_2^* W  \left( \begin{smallmatrix}   &  & 1 \\   & 1  & \\ 1 &  &    \end{smallmatrix} \right)Z^* = V_1V_1^* ,
 $$
and hence that $TV_2V_2^* T^* = V_2V_2^*$. It follows then from
Lemma \ref{B69} that $T$, and hence also $ZW^*$ is connected to
$1$ in the unitary group of $\beta_0(A)'\cap M(B)$.
\end{proof}
\end{lemma}

\begin{lemma}\label{B94}
Let $\varphi, \psi, \lambda : A \to Q(B)$ be extensions of $A$ by
$B$ that equal $\alpha$ on $D$. There are then unitaries $v,w \in
\alpha(A)'\cap Q(B)$, connected to $1$ in the unitary group of $
\alpha(A)'\cap Q(B),$ such that $\Ad v \circ \left( \varphi
+_{\alpha} \psi\right) = \psi +_{\alpha} \varphi$ and $\Ad w \circ
\left( \left( \varphi +_{\alpha} \psi\right) +_{\alpha}
\lambda\right) = \varphi +_{\alpha} \left( \psi +_{\alpha}
\lambda\right)$.
\begin{proof}
By the essential uniqueness of absorbing $*$-homomorphisms there
is a unitary $U \in M(B)$ such that $\Ad U \circ \beta_0(a) -
\alpha_0(a) \in B$. Set $t = w_1q_B(UV_1U^*)^* +
w_2q_B(UV_2U^*)^*$, and let $S,T$ be the unitaries from Lemma
\ref{B73}. Then $v = tq_B(USU^*)t^*$ and $w = tq_B(UTU^*)t^*$ have
the stated properties by Lemma \ref{B73}.
\end{proof}
\end{lemma}

It follows from Lemma \ref{B94} that (\ref{B92}) gives the set of
unitary equivalence classes of extensions of $A$ by $B$ which
agree with $\alpha$ on $J$ the structure of an abelian semi-group.
Furthermore, it follows from the proof of Lemma \ref{B94} that any
other choice of absorbing $*$-homomorphism instead of $\alpha$
would result in an isomorphic semi-group. To obtain a neutral
element we declare that two extensions $\varphi, \psi : A \to
Q(B)$  that agree with $\alpha$ on $J$ are \emph{stably
equivalent} when $\varphi +_{\alpha} \alpha$ is unitarily
equivalent to $\psi +_{\alpha} \alpha$. That stable equivalence
\emph{is} an equivalence relation follows from Lemma \ref{B94} and
the observation that $\alpha +_{\alpha} \alpha = \alpha$. The
formula (\ref{B92}) gives a well-defined composition in the set of
stable equivalence classes of extensions that agree with $\alpha$
on $J$, giving us an abelian semi-group with a neutral element (or
$0$) represented by $\alpha$, and we denote this semi-group by
$\Ext_{J,\alpha}(A,B)$. The group of invertible elements in
$\Ext_{J,\alpha}(A,B)$ will be denoted by
 $$
\Ext_{J,\alpha}^{-1}(A,B) .
 $$
It is clear from the construction, and can be seen from the
essential uniqueness of an absorbing $*$-homomorphism, that any
other choice of an absorbing $*$-homomorphism $A \to M(B)$ will
give rise to an isomorphic group. However, at this point it would
seem as if the stable equivalence of two given extensions of $A$
by $B$, which both agree with $\alpha$ on $J$, depends on the
particular choice of isometries from $\alpha(A)'\cap Q(B)$ used to
define the addition $+_{\alpha}$. The next lemma shows that this
is not the case because the addition (\ref{B92}) is independent of
the $w_i$'s up to conjugation by a unitary from $ \alpha(A)'\cap
Q(B)$.

Let $n \in \mathbb N$. To simplify the notation, we denote by $1_n
\otimes \alpha_0:A \to M_n\left(M(B)\right)$ the $*$-homomorphism
given by
 $$
\left(1_n \otimes \alpha_0\right)(a) =  \left( \begin{smallmatrix} \alpha_0(a) & & &  \\  & \alpha_0(a)  &  & \\   & & \ddots & \\ & & &  \alpha_0(a) \end{smallmatrix}\right).
 $$
Set $1_n \otimes \alpha = \left(\id_{M_n(\mathbb C)} \otimes
q_B\right) \circ \left(1_n \otimes \alpha_0\right)$.

\begin{lemma}\label{ideal2}
Let $\varphi : A \to Q(B)$ be an extension of $A$ by $B$ which is
equal to $\alpha$ on $J$, and assume that $v \in \alpha(A)'\cap
Q(B)$ is an isometry such that $vv^* \alpha(a) = \alpha(a)$ for
all $a \in A$. It follows that $\Ad v \circ \varphi$ is stably
equivalent to $\varphi$.
\begin{proof} Note that
\begin{equation*}\label{ideal3}
\left( \begin{matrix} v & 1 - vv^* \\ 0  & v^* \end{matrix}\right) \left( \begin{matrix} \varphi(a)  & 0 \\ 0  & \alpha(a) \end{matrix}\right)\left( \begin{matrix} v^* &  0 \\ 1 - vv^* & v \end{matrix}\right) =  \left( \begin{matrix} v\varphi(a)v^*  & 0 \\ 0  & \alpha(a) \end{matrix}\right)
\end{equation*}
and that $\left( \begin{matrix} v & 1 - vv^* \\ 0  & v^*
\end{matrix}\right)$ is a unitary in $M_2(Q(B)) \cap \left(1_2
\otimes \alpha\right)(A)'$. It follows that there is a unitary $u
\in Q(B) \cap \alpha(A)'$ such that $\Ad w \circ \left( \varphi
+_{\alpha} \alpha\right) = (\Ad v \circ \varphi) +_{\alpha}
\alpha$. Since $\left( \begin{smallmatrix} u &  \\   & u^*
\end{smallmatrix}\right)$ is connected to $1$ in the unitary group
of $M_2(Q(B)) \cap \left(1_2 \otimes \alpha\right)(A)'$, we deduce
that $\Ad v \circ \varphi +_{\alpha} \alpha  +_{\alpha} \alpha$ is
unitarily equivalent to $\varphi  +_{\alpha} \alpha  +_{\alpha}
\alpha$.
\end{proof}
\end{lemma}

\section{A six-term exact sequence.}\label{sixterms}

We will now assume that there is an absorbing $*$-homomorphism
$\alpha_0 : A \to M(B)$ such that $\alpha_0|_J : J \to M(B)$ is
also absorbing. This condition is known to be automatically
fullfilled in the following cases:
\begin{enumerate}
\item[i)] $B$ is nuclear, or
\item[ii)] $J$ is nuclear, or
\item[iii)]$J$ is a hereditary $C^*$-subalgebra of $A$; in particular, when $J$ is an ideal, or
\item[iv)] there is a surjective conditional expectation $P : A \to J$.
\end{enumerate}
i) follows from Kasparov's work, \cite{K2}, and ii)-iv) all follow from Lemma 2.1 and Lemma 2.2 of \cite{Th3}. In general the existence of $\alpha_0$ fails, cf. the last section in \cite{Th3}.

\smallskip

Fix an absorbing $*$-homomorphism $\alpha_0 : A \to M(B)$ such
that $\alpha_0|_J : J \to M(B)$ is also absorbing, and set $\alpha
= q_B \circ \alpha_0$ as before. Set
 $$
\mathcal D_{\alpha}(J) = \left\{ m \in M(B): \ m \alpha_0(j) -
\alpha_0(j)m \in B \ \forall j \in J \right\}
 $$
and
 $$
\mathcal X_{\alpha}(J) = \left\{ m \in D_{\alpha}(J): \ m \alpha_0(j) \in B \ \forall j \in J \right\} .
 $$
It was shown in \cite{Th1} that there is a natural isomorphism
$K_1\left(\mathcal D_{\alpha}(J)/\mathcal X_{\alpha}(J)\right)
\simeq KK(J,B)$, and then in Lemma 3.1 of \cite{Th2} that
$K_*\left(\mathcal X_{\alpha}(J)\right) = 0$, so that we have a
natural isomorphism
\begin{equation}\label{B32}
K_1\left(\mathcal D_{\alpha}(J) \right) \simeq KK(J,B) .
\end{equation}
Similarly, we set $\mathcal D_{\alpha}(A) = \left\{ m \in M(B): \ m \alpha_0(a) -  \alpha_0(a)m \in B \ \forall a \in A \right\}$, and get a natural isomorphism
\begin{equation}\label{B51}
K_1\left(\mathcal D_{\alpha}(A) \right) \simeq KK(A,B) .
\end{equation}
As above $i : J \to A$ will denote the inclusion, and we denote
also by $i$ the induced inclusion $i: \mathcal D_{\alpha}(A) \to
\mathcal D_{\alpha}(J)$. Note that the diagram
\begin{equation}\label{B34}
\begin{xymatrix}{
K_1\left( \mathcal D_{\alpha}(A)   \right) \ar[d] \ar[r]^-{i_*} & K_1\left(\mathcal D_{\alpha}(J) \right) \ar[d] \\
KK(A,B) \ar[r]^-{i^*} & KK(J,B) }
\end{xymatrix}
\end{equation}
is commutative when the vertical arrows are the isomorphisms (\ref{B32}) and (\ref{B51}).

Let $v$ be a unitary in $M_n\left(\mathcal D_{\alpha}(J)\right)$.
Let $\Theta_n : M_n(Q(B)) \to Q(B)$ be the isomorphism from Lemma
\ref{B91}. Then $\Theta_n \circ \left(\id_{M_n(\mathbb C)} \otimes
q_B\right) \circ \Ad v \circ \left(1_n \otimes \alpha_0\right) : A
\to M_n(Q(B))$ is an extension $e(v) : A \to Q(B)$ of $A$ by $B$
which is equal to $\alpha$ on $J$. Since
 $$
\Theta_{2n} \circ \left(\id_{M_{2n}(\mathbb C)} \otimes q_B\right)\left( \begin{smallmatrix} v  & \\ & v^*\end{smallmatrix} \right)
 $$
is connected to $1$ in the unitary group of $\alpha(J)'\cap Q(B)$,
we see that $e(v) +_{\alpha} e(v^*)$ is stably equivalent to
$\alpha$, as an extension of $A$ by $B$ which is equal to $\alpha$
on $J$, proving that $e(v)$ is invertible, hence it represents an
element in $\Ext^{-1}_{J,\alpha}(A,B)$. When $v_t, t \in [0,1]$,
is a norm-continuous path of unitaries in $M_n\left(\mathcal
D_{\alpha}(J)\right)$ there is a partition $0 =t_0 < t_1< t_2 <
\dots < t_N =1$ of $[0,1]$ such that $\left(\id_{M_n(\mathbb C)}
\otimes q_B\right)\left(v_{t_i}v_{t_{i+1}}^*\right)$ is in the
connected component of $1$ in the unitary group of $M_n\left(
\alpha(J)' \cap Q(B)\right)$. Hence $e\left(v_0\right) =
e\left(v_1\right)$ in $\Ext^{-1}_{J,\alpha}(A,B)$. It is then
straightforward to check that the construction gives us a group
homomorphism
\begin{equation}\label{rom1}
\partial : K_1\left(\mathcal D_{\alpha}(J) \right) \to \Ext^{-1}_{J,\alpha}(A,B) .\end{equation}

\begin{lemma}\label{B35} The sequence
\begin{equation*}
\begin{xymatrix}{
 K_1\left(\mathcal D_{\alpha}(J)\right) \ar[r]^-{\partial}  & \Ext^{-1}_{J,\alpha}(A,B) \ar[r] & \Ext^{-1}(A,B)  \ar[d]^-{i^*} \\
K_1\left(\mathcal D_{\alpha}(A)\right) \ar[u]^-{i_*} &  &     \Ext^{-1}(J,B)}
\end{xymatrix}
\end{equation*}
is exact.
\begin{proof}
Exactness at $K_1\left(\mathcal D_{\alpha}(J)\right) $: If $v$ is
a unitary in $M_n\left( \mathcal D_{\alpha}(A)\right)$, we find
that $\left(\id_{M_n(\mathbb C)} \otimes q_B\right) \circ \Ad v
\circ \left(1_n \otimes \alpha_0\right) = 1_n \otimes \alpha$ so
that $e(v)=[\Theta_n\circ(1_n\otimes\alpha)]=[\alpha] = 0$ in
$\Ext^{-1}_{J,\alpha}(A,B)$. To show that $\ker \partial \subseteq
\im i^*$, let $v \in \mathcal D_{\alpha}(J)$ be a unitary such
that $\partial[v] = 0$. Then $ \Ad q_B(v) \circ \alpha +_{\alpha}
\alpha$ is unitarily equivalent to $\alpha  +_{\alpha} \alpha$,
which means that there is a unitary $S$ connected to $1$ in the
unitary group of $M_2\left(\alpha(J)'\cap Q(B)\right)$ such that
\begin{equation}\label{B6}
\Ad \left[S\left( \begin{smallmatrix} q_B(v) & {} \\ {} & 1 \end{smallmatrix} \right)\right]  \left(\begin{smallmatrix}  \alpha(a) & {} \\ {} & \alpha(a) \end{smallmatrix} \right )    =  \left ( \begin{smallmatrix}\alpha(a) & {} \\ {} & \alpha(a) \end{smallmatrix} \right ).
\end{equation}
Since $q_B : \mathcal D_{\alpha}(J) \to \alpha(J)'\cap Q(B)$ is
surjective, there is a unitary $S_0$ connected to $1$ in the
unitary group of $M_2\left(\mathcal D_{\alpha}(J)\right)$ such
that $\id_{M_2(\mathbb C)} \otimes q_B(S_0) = S$. Then $[v] =
\left[S_0\left( \begin{smallmatrix} v & {} \\ {} & 1
\end{smallmatrix} \right)\right]$ in $K_1\left(\mathcal
D_{\alpha}(J)\right)$. It follows from (\ref{B6}) that $S_0\left(
\begin{smallmatrix} v & {} \\ {} & 1
\end{smallmatrix} \right)$ belongs to the commutant of $(1_2\otimes\alpha)(A)$,
hence $\left[S_0\left( \begin{smallmatrix} v & {} \\ {} & 1
\end{smallmatrix} \right)\right] \in i^*\left(K_1\left(\mathcal
D_{\alpha}(A)\right)\right)$. The same argument works when $v$ is
a unitary in $M_n\left(\mathcal D_{\alpha}(J)\right)$ for some $n
\geq 2$.

\smallskip

Exactness at $ \Ext^{-1}_{J,\alpha}(A,B)$: Let $v$ be a unitary in
$\mathcal D_{\alpha}(J)$. Then $(\Ad q_B(v) \circ  \alpha)$ is a
split extension of $A$ by $B$, proving that $[v] \in K_1\left(
\mathcal D_{\alpha}(J)\right)$ goes to zero under the composition
$ K_1\left(\mathcal D_{\alpha}(J)\right) \to
\Ext^{-1}_{J,\alpha}(A,B) \to \Ext^{-1}(A,B)$. The same argument
works when $v \in M_n \left(\mathcal D_{\alpha}(J)\right)$, so we
see that the composition is zero. Let $\varphi : A \to Q(B)$ be an
extension of $A$ by $B$ which is equal to $\alpha$ over $J$, and
assume that $[\varphi] = 0$ in $\Ext^{-1}(A,B)$. This means that
there is a unitary $T \in M(M_2(B))$ such that
\begin{equation}\label{hjt1}
\Ad \left( \id_{M_2(\mathbb C)} \otimes q_B\right)(T) \circ  \left ( \begin{smallmatrix} \varphi & {}  \\ {} &  \alpha   \end{smallmatrix} \right ) =  \left ( \begin{smallmatrix}  \alpha & {}  \\ {} &  \alpha  \end{smallmatrix} \right ) .
\end{equation}
Since $\varphi$ is equal to $\alpha$ over $J$ this implies that $T
\in M_2\left(\mathcal D_{\alpha}(J)\right)$ and we see from
(\ref{hjt1}) that $[\varphi] = \partial[T^*]$. Note that we did
not assume that $\varphi$ represented an \emph{invertible} element
in $\Ext_{J,\alpha}(A,B)$, so besides establishing the exactness
at $ \Ext^{-1}_{J,\alpha}(A,B)$ the argument also shows that every
element of $\Ext_{J,\alpha}(A,B)$ which goes to $0$ in
$\Ext^{-1}(A,B)$ comes from $K_1\left(\mathcal
D_{\alpha}(J)\right)$, and hence \emph{is} invertible in
$\Ext_{J,\alpha}(A,B)$. This point will be used shortly.

\smallskip

Exactness at $ \Ext^{-1}(A,B)$: It is obvious that $i^*$ kills the
image of $\Ext^{-1}_{J}(A,B)$ in $\Ext^{-1}(A,B)$, so consider an
invertible extension $\varphi : A \to Q(B)$ such that $[\varphi
\circ i] = 0$ in $\Ext^{-1}(J, B)$. This means that there is a
unitary $T \in M(B)$ such that
 $$
\Ad q_B(T)\left ( \varphi \oplus \alpha\right)(j) = \alpha(j)
 $$
for any $j\in J$, i.e $\Ad q_B(T)\left ( \varphi \oplus
\alpha\right)$ equals $\alpha$ on $J$. Since $[\varphi] =\left[\Ad
q_B(T)\left ( \varphi \oplus \alpha\right)\right]$ in
$\Ext^{-1}(A,B)$, this completes the proof, provided we can show
that $\Ad q_B(T)\left ( \varphi \oplus \alpha\right)$ represents
an invertible element of $\Ext_{J,\alpha}(A,B)$. To this end, let
$\psi : A \to Q(B)$ be an extension of $A$ by $B$ which represents
the inverse of $\varphi$ in $\Ext^{-1}(A,B)$. Then $i^*[\psi] = 0$
in $\Ext^{-1}(J,B)$ so we deduce as in case of $\varphi$ that
there is a unitary $T'\in M(B)$ such that $\Ad q_B(T')\left ( \psi
\oplus \alpha\right)(j) = \alpha(j)$ for all $j \in J$. Thus $\Ad
q_B(T')\left ( \psi \oplus \alpha\right)$ and $\Ad q_B(T)\left (
\varphi \oplus \alpha\right)$ both represent elements of
$\Ext_{J,\alpha}(A,B)$. Since the sum
 $$
\Ad q_B(T')\left ( \psi \oplus \alpha\right)  +_{\alpha} \Ad q_B(T)\left ( \varphi \oplus \alpha\right)
 $$
represents $0$ in $\Ext^{-1}(A,B)$, it follows from the argument
that proved exactness at $ \Ext^{-1}_{J,\alpha}(A,B)$ that $\left[
\Ad q_B(T)\left ( \varphi \oplus \alpha\right)\right]$ is
invertible in $\Ext_{J,\alpha}(A,B)$.
\end{proof}
\end{lemma}

The proof of Lemma \ref{B35} has the following corollary:

\begin{lemma}\label{C35}
Let $\varphi$ be an extension of $A$ by $B$ which equals $\alpha$
on $J$. Assume that $\varphi$ is invertible in $\Ext(A,B)$. It
follows that $\varphi$ is invertible in $\Ext_{J,\alpha}(A,B)$.
\begin{proof}
The image of $[\varphi]$ in $\Ext^{-1}(A,B)$ is killed by $i^*:
\Ext^{-1}(A,B) \to \Ext^{-1}(J,B)$ and the proof Lemma \ref{B35},
more precisely the proof of exactness at $\Ext^{-1}(A,B)$, shows
that $[\varphi]$ is invertible in $\Ext_{J,\alpha}(A,B)$.
\end{proof}
\end{lemma}

In particular, $\Ext_{J,\alpha}(A,B)$ is a group when $\Ext(A,B)$ is.

\smallskip

Consider now the suspension $SB = C_0(0,1) \otimes B$ of $B$. If
we combine (\ref{B34}) with Lemma \ref{B35} and use the natural
identification $KK(-,B) = \Ext^{-1}(-,SB)$ we get the exact
sequence
\begin{equation}\label{B36}
\begin{xymatrix}{
 \Ext^{-1}(J,SB) \ar[r]^-{\partial_0}  & \Ext^{-1}_{J,\alpha}(A,B) \ar[r] & \Ext^{-1}(A,B)  \ar[d]^-{i^*} \\
\Ext^{-1}\left(A,SB\right) \ar[u]^-{i^*} &  &     \Ext^{-1}(J,B)}
\end{xymatrix}
\end{equation}
where $\partial_0$ is the composition of $\partial :
K_1\left(D_{\alpha}(J)\right) \to \Ext^{-1}_J(A,B)$ with the
isomorphism $\Ext^{-1}(J,SB) \to  K_1\left( D_{\alpha}(J)\right)$
coming from (\ref{B32}).

Let now $\beta_0 : A \to M(SB)$ be an absorbing $*$-homomorphism
such that also $\beta_0|_J :J \to M(SB)$ is absorbing. The
existence of such a $\beta_0$ does not require any additional
asumptions because $\beta_0$ can be constructed from $\alpha_0$ by
use of Lemma 3.2 of \cite{Th2}. Hence there is also a map
$\partial': K_1\left(\mathcal D_{\beta}(J)\right) \to
\Ext_{J,\beta}^{-1}(A,SB)$, defined in the same way as $\partial$,
but with $SB$ in place of $B$. This leads to the following version
of (\ref{B36}):
\begin{equation}\label{B36'}
\begin{xymatrix}{
\Ext^{-1}\left(J,SB\right)  &   &\Ext^{-1}\left(A,B\right) \ar[d]^-{i^*} \\
\Ext^{-1}(A,SB) \ar[u]^-{i^*}  & \Ext^{-1}_{J,\beta}(A,SB) \ar[l] &  \Ext^{-1}(J,B) \ar[l]^-{\partial_1}, }
\end{xymatrix}
\end{equation}
where $\partial_1$ is the composition of $\partial':
K_1\left(\mathcal D_{\beta}(J)\right) \to
\Ext^{-1}_{J,\beta}(A,SB)$ with the isomorphism $\Ext^{-1}(J,B)=
KK(J,SB) \to  K_1\left(\mathcal D_{\beta}(J)\right)$. By combining
(\ref{B36}) and (\ref{B36'}) we get

\begin{thm}\label{cor1} The sequence
\begin{equation*}
\begin{xymatrix}{
 \Ext^{-1}_{J,\alpha}(A,B)  \ar[r]  &   \Ext^{-1}(A,B)  \ar[r]^-{i^*} & \Ext^{-1}(J,B)  \ar[d]^-{\partial_1} \\
\Ext^{-1}(J,SB)  \ar[u]^-{\partial_0} &   \Ext^{-1}\left(A,SB\right) \ar[l]^-{i^*} &  \Ext^{-1}_{J,\beta}(A,SB) \ar[l] }
\end{xymatrix}
\end{equation*}
is exact.
\end{thm}

\section{Other realizations of the relative extension group.}\label{other}

As above we assume that there is an absorbing $*$-homomorphism
$\alpha_0: A \to M(B)$ such that $\alpha_0|_J : J \to M(B)$ is
also absorbing. By an \emph{extension of $A$ by $B$ which splits
over $J$} we mean a pair $(\varphi, \varphi_0)$ where $\varphi : A
\to Q(B)$ is an extension of $A$ by $B$, and $\varphi_0 : J \to
M(B)$ is a $*$-homomorphism such that $q_B \circ \varphi_0 =
\varphi|_J$. We say that $(\varphi, \varphi_0)$ is
\emph{invertible} when $\varphi$ is an invertible extension of $A$
by $B$, i.e. when there is another extension $\psi$ of $A$ by $B$
with the property that $\varphi \oplus \psi$ is a split extension
(of $A$ by $B$). Two invertible extensions, $(\varphi,\varphi_0)$
and $(\psi,\psi_0)$, of $A$ by $B$ which split over $J$ are
\emph{homotopic in norm} when there is a path
$\left(\Phi^t,\Phi_0^t\right), t \in [0,1]$, of invertible
extensions of $A$ by $B$ which split over $J$ such that $[0,1] \ni
t \mapsto \Phi^t(a)$ is norm-continuous for all $a \in A$ and
$[0,1] \ni t \mapsto \Phi_0^t(j)$ is norm-continuous for all $j
\in J$, $\left(\Phi^0,\Phi_0^0\right) = (\varphi,\varphi_0)$ and
$\left(\Phi^1,\Phi_0^1\right) = (\psi,\psi_0)$. We say that
$(\varphi,\varphi_0)$ and $(\psi,\psi_0)$ are \emph{stably
homotopic in norm} when there is a $*$-homomorphism $\pi : A \to
M(B)$ such that $\left( \varphi \oplus q_B \circ \pi, \varphi_0
\oplus \pi|_J\right)$ and $\left( \psi \oplus q_B \circ \pi,
\psi_0 \oplus \pi|_J\right)$ are homotopic in norm. We denote by
$\Ext_{J}^{-1}(A,B)$ the abelian semi-group of stable homotopy
classes of invertible extensions of $A$ by $B$ which split over
$J$. As we shall see shortly, $\Ext_{J}^{-1}(A,B)$ is actually a
group.

Choose a sequence $W_1,W_2,W_3, \dots $ of isometries in $M(B)$
such that $W_i^*W_j =0$ when $i \neq j$, and $\sum_{i=1}^{\infty}
W_iW_i^* = 1$, with convergence in the strict topology. Set
$\beta_0(a) = \sum_{i=2}^{\infty} W_i \alpha_0(a)W_i^*$, and note
that $\beta_0: A \to M(B)$ and $\beta_0|_J : J \to M(B)$ are both
absorbing. We shall work with $\beta_0$ instead of $\alpha_0$. The
point is that unlike $\alpha_0$, the absorbing $*$-homomorphisms
$\beta_0$ and $\beta_0|_J$ are both guaranteed to be saturated in
the sense of \cite{Th2}. Recall that a $*$-homomorphism
$\varphi_0:A\to M(B)$ is {\it saturated} if it is unitarily
equivalent to $(\varphi_0)_\infty\oplus(0)_\infty$, where
$(\varphi_0)_\infty=\varphi_0\oplus\varphi_0\oplus\ldots$.

Let $\varphi : A \to Q(B)$ be an extension of $A$ by $B$ which
equals $\beta$ on $J$. Then $(\varphi, \beta_0)$ is an extension
of $A$ by $B$ which splits over $J$, and it is straightforward to
see that the recipe $[\varphi] \to [\varphi, \beta_0]$ is a group
homomorphism
 \begin{equation}\label{iso41}
\Ext^{-1}_{J,\beta}\left( A,B\right) \to \Ext_J^{-1}(A,B).
\end{equation}
The aim is to show that (\ref{iso41}) is an isomorphism.

Set $IB = C[0,1] \otimes B$, and let $\ev_t : IB \to B$ be the
$*$-homomorphism given by evaluation at $t \in [0,1]$. Then $e_t$
extends to a $*$-homomorphism $\overline{\ev_t} : M(IB)\to M(B)$
and induces in turn a $*$-homomorphism $\widehat{\ev_t} : Q(IB)
\to Q(B)$. Let $\gamma_0 : A \to M(IB)$ be the $*$-homomorphism
such that $\left(\gamma_0(a)f\right)(t) = \beta_0(a)f(t), t \in
[0,1], f \in IB$. Since $\beta_0$ is saturated it follows from
Lemma 2.3 of \cite{Th2} that $\gamma_0$ is absorbing. Set $\gamma
= q_{IB} \circ \gamma_0$, and note that we have, for any $t \in
[0,1]$, a homomorphism
\begin{equation*}\label{homtop}
{e_t}_* : \Ext_{J,\gamma}^{-1}(A,IB) \to \Ext_{J,\beta}^{-1}(A,B)
\end{equation*}
defined such that ${e_t}_*[\varphi] = \left[\widehat{e_t} \circ
\varphi\right]$ when $\varphi : A \to Q(IB)$ is an extension of
$A$ by $IB$ which equals $\gamma$ on $J$.

\begin{lemma}\label{homtop2}
The homomorphisms ${e_t}_*, t\in [0,1]$, are all the same group
isomorphism.
\begin{proof}
Define $c : B \to IB$ such that $c(b)(t) = b$, and note that $c$
induces $*$-homomorphisms $\overline{c} : M(B) \to M(IB)$ and
$\widehat{c}: Q(B) \to Q(IB)$. Since $\widehat{c} \circ \beta =
\gamma$ there is a homomorphism $c_* : \Ext_{J,\beta}^{-1}(A,B)
\to \Ext_{J,\gamma}^{-1}(A,IB)$ such that $c_*[\psi] =
[\widehat{c} \circ \psi]$. Since ${e_t}_* \circ c_*$ is the
identity on $\Ext_{J,\beta}^{-1}(A,B)$ for all $t \in [0,1]$, it
suffices to show that $c_*$ is an isomorphism. This follows from
Theorem \ref{cor1} by an obvious application of the five-lemma.
\end{proof}
\end{lemma}

\begin{lemma}\label{iso6}
Let $A_1$ and $B_1$ be separable $*$-algebras, $B_1$ stable. Let
$\varphi, \psi : A_1 \to M(B_1)$ be $*$-homomorphisms and $W_t, t
\in [1,\infty)$, a continuous path of unitaries in $M(B_1)$ such
that
\begin{enumerate}
\item[i)] $W_t \varphi(a)W_t^* - \psi(a) \in A_1, \ t \in [1,\infty), a \in A_1$,
\item[ii)] $\lim_{t \to \infty} W_t\varphi(a)W_t^* = \psi(a), \ a \in A_1$.
\end{enumerate}
Then $\left[\Ad W_1 \circ \varphi, \psi\right] = 0$ in $KK(A_1,B_1)$.
\begin{proof}
The lemma and its proof are essentially identical to Lemma 3.1 of
\cite{DE2}. Note, however, that one of the crucial assumptions has
mysteriously disappeared in the lemma in \cite{DE2}.
\end{proof}
\end{lemma}

\begin{prop}\label{iso42}
The map (\ref{iso41}) is an isomorphism. In particular,
$\Ext_J^{-1}(A,B)$ is a group.
\begin{proof}
Surjectivity: Let $(\varphi,\varphi_0)$ be an invertible extension
of $A$ by $B$ which splits over $J$. Then $\varphi_0 \oplus
\beta_0$ is approximately unitarily equivalent to $\beta_0$
because $\beta_0$ is absorbing, i.e. there exists a sequence of
unitaries $U_n\in M(B)$, $n\in\mathbb N$, such that $\Ad
U_n\circ(\varphi_0\oplus\beta_0)(j)-\beta_0(j)\in B$ for all $n$,
and $\lim_{n\to\infty}\Ad
U_n\circ(\varphi_0\oplus\beta_0)(j)-\beta_0(j)=0$ for any $j\in
J$. It follows then from Lemma
2.4 of \cite{DE2} that $\varphi_0 \oplus
\left(\beta_0\right)_{\infty}$ is asymptotically unitarily
equivalent to $\left(\beta_0\right)_{\infty}$, i.e. there exists a
norm-continuous path $V_t, t \in [1,\infty)$, of unitaries in
$M(B)$ such that $\Ad
V_t\circ(\varphi_0\oplus\beta_0)(j)-\beta_0(j)\in B$ for all $t$,
and $\lim_{t\to\infty}\Ad
V_t\circ(\varphi_0\oplus\beta_0)(j)-\beta_0(j)=0$ for any $j\in
J$. Since $\beta_0$ is unitarily equivalent to $
\left(\beta_0\right)_{\infty}$ because $\beta_0$ is saturated, we
conclude that there is a norm-continuous path $W_t, t \in
[1,\infty)$, of unitaries in $M(B)$ such that $W_t\left(\varphi_0
\oplus \beta_0\right)(j)W_t^* - \beta_0(j) \in B$ for all $t$, and
$\lim_{t \to \infty} W_t\left(\varphi_0 \oplus
\beta_0\right)(j)W_t^* - \beta_0(j) = 0$ for any $j \in J$. Since
the unitary group of $M(B)$ is connected in norm, it holds that
$[\varphi,\varphi_0] = \left[\Ad q_B(W_1) \circ \left(\varphi
\oplus \beta\right), \Ad W_1 \circ \left(\varphi_0 \oplus
\beta_0\right)\right]$ in $ \Ext_{J,s}(A,B)$. Set
 $$
\Psi^t = \begin{cases} \Ad W_{\frac{1}{t}} \circ \left(\varphi_0 \oplus \beta_0\right), & t \in ]0,1] \\ \beta_0 , & t = 0  . \end{cases}
 $$
Then $\left(\Ad q_B(W_1) \circ \left(\varphi \oplus \beta\right),
\Psi^t\right), t \in [0,1]$, is a homotopy in norm showing that
 $$
\left[\Ad q_B(W_1) \circ \left(\varphi \oplus \beta\right), \Ad W_1 \circ \left(\varphi_0 \oplus \beta_0\right)\right] = \left[\Ad q_B(W_1) \circ \left(\varphi \oplus \beta\right), \beta_0\right].
 $$
Note that $\Ad q_B(W_1) \circ \left(\varphi \oplus \beta\right)$
is equal to $\beta$ on $J$, and is invertible in
$\Ext_{J,\beta}(A,B)$ by the proof of Lemma \ref{B35} since
$\varphi$ is an invertible extension. Thus $\Ad q_B(W_1) \circ
\left(\varphi \oplus \beta\right)$ represents an element of
$\Ext^{-1}_{J,\beta}(A,B)$, and we conclude that (\ref{iso41}) is
surjective. In particular $\Ext^{-1}_{J}(A,B)$ is a group.

Injectivity: Let $\varphi,\psi$ be extensions of $A$ by $B$ which
both equal $\beta$ on $J$. Assume that $[\varphi,\beta_0] =
[\psi,\beta_0]$ in $\Ext_{J}^{-1}(A,B)$. There is then a
$*$-homomorphism $\pi : A \to M(B)$ such that $\left( \varphi
\oplus q_B \circ \pi, \beta_0 \oplus \pi\right)$ is homotopic in
norm to $\left( \psi \oplus q_B \circ \pi, \beta_0 \oplus
\pi\right)$. As in the proof of surjectivity we can find a
norm-continuous path $W_t, t \in [1,\infty)$, of unitaries in
$M(B)$ such that $W_t\left(\pi\oplus \beta_0\right)(a)W_t^* -
\beta_0(a) \in B$ for all $t$, and $\lim_{t \to \infty}
W_t\left(\pi \oplus \beta_0\right)(a)W_t^* - \beta_0(a) = 0$ for
all $a \in A$. Set $S_t = 1 \oplus W_t$ and note that this gives a
homotopy in norm between $( \varphi \oplus q_B\circ \pi \oplus
\beta,  \beta_0 \oplus \pi \oplus \beta)$ and $\left(\varphi
\oplus \beta,  \beta_0 \oplus \beta_0\right)$. Similarly, there is
homotopy in norm between $( \psi \oplus q_B\circ \pi \oplus \beta,
\beta_0 \oplus \pi \oplus \beta)$ and $\left(\psi \oplus \beta,
\beta_0 \oplus \beta_0\right)$. It follows that there is a
homotopy in norm, $\left(\Psi^t,\Psi_0^t\right)$, between
$\left(\varphi +_{\beta} \beta,  \beta_0 +_{\beta} \beta_0\right)$
and $\left(\psi +_{\beta} \beta,  \beta_0 +_{\beta}
\beta_0\right)$. This homotopy defines in an obvious way an
extension $\Phi : A \to Q(IB)$ and a $*$-homomorphism $\Phi_0 : J
\to M(IB)$ such that $q_{IB} \circ \Phi_0 = \Phi|_J$,
$(\widehat{e_0} \circ \Phi, \overline{e_0} \circ \Phi_0) =
\left(\varphi +_{\beta} \beta,  \beta_0 +_{\beta} \beta_0\right)$
and $(\widehat{e_1} \circ \Phi, \overline{e_1} \circ \Phi_0) =
\left(\psi +_{\beta} \beta,  \beta_0 +_{\beta} \beta_0\right)$ .
Note that $\Phi$ is invertible since each $\Psi^t$ is. Let
$\gamma_0 : A \to M(IB)$ be the $*$-homomorphism such that
$\left(\gamma_0(a)f\right)(t) = \beta_0(a)f(t), t \in [0,1], f \in
IB$. As in the proof of surjectivity, it follows from \cite{DE2}
that there is a norm-continuous path $W_t, t \in [1,\infty)$, of
unitaries in $M(IB)$ such that $W_t\left(\Phi_0 \oplus
\gamma_0\right)(j)W_t^* - \gamma_0(j) \in IB$ for all $t$, and
$\lim_{t \to \infty} W_t\left(\Phi_0 \oplus
\gamma_0\right)(j)W_t^* - \gamma_0(j) = 0$ for any $j \in J$.
There is therefore also  a norm-continuous path $V_t, t \in
[1,\infty)$, of unitaries in $M(IB)$ such that $V_t\left(\Phi_0
{+}_{\gamma} \gamma_0\right)(j)V_t^* - \left(\left(\gamma_0
+_{\gamma} \gamma_0 \right) +_{\gamma} \gamma_0\right)(j) \in IB$
for all $t$, and $\lim_{t \to \infty} V_t\left(\Phi_0 +_{\gamma}
\gamma_0\right)(j)V_t^* - \left(\left(\gamma_0 +_{\gamma} \gamma_0
\right) +_{\gamma} \gamma_0\right)(j) = 0$ for any $j \in J$. It
follows from Lemma \ref{iso6} that
\begin{equation}\label{iso7}
\left[\Ad V_1 \circ \left(\Phi_0 +_{\gamma} \gamma_0\right), \left(\gamma_0 +_{\gamma} \gamma_0 \right) +_{\gamma} \gamma_0  \right] = 0
\end{equation}
in $KK(J  ,IB )$. Set $S_i = \overline{e_i}(V_1), i = 0,1$, and
note that it follows from (\ref{iso7}) that
\begin{equation}\label{iso8}
 \left[\Ad S_0 \circ \left( \left( \beta_0 +_{\beta} \beta_0\right)  +_{\beta} \beta_0\right), \beta_0 +_{\beta} \beta_0\right)  +_{\beta} \beta_0  ] =  0
\end{equation}
and
\begin{equation}\label{iso9}
 \left[\Ad S_1 \circ \left( \left( \beta_0 +_{\beta} \beta_0\right)  +_{\beta} \beta_0\right), \left(\beta_0 +_{\beta} \beta_0\right)  +_{\beta} \beta_0  \right] = 0
\end{equation}
in $KK(J,B)$. Let $W \in M(B)$ be a unitary such that $\Ad W \circ
\beta_0 = \left(\beta_0 +_{\beta} \beta_0\right)$. It follows then
from (\ref{iso8}) and (\ref{iso9}) that
\begin{equation}\label{iso8'}
\left[\Ad W^*S_0W \circ \beta_0, \beta_0\right] =\left[\Ad W^*S_1W \circ \beta_0, \beta_0\right]  = 0
\end{equation}
in $KK(J,B)$. Note that $\left[\Ad W^*S_0W \circ \beta_0,
\beta_0\right] \in KK(J,B)$ is the image of the class of the
unitary $W^*S_0W$ under the isomorphism (\ref{B32}). Thus
(\ref{iso8'}) implies that $\left[q_B(W^*S_0W)\right] = 0$ in
$K_1\left( \beta\left(J\right)^{\prime} \cap Q(B)\right)$.
Similarly, it also implies that  $\left[q_B(W^*S_1W)\right] = 0$
in $K_1\left( \beta\left(J\right)^{\prime} \cap Q(B)\right)$. It
follows therefore that
 $$
[\varphi] = \left[\Ad q_B(W^*S_0W) \circ  \left( \widehat{e_0} \circ \Phi  +_{\beta} \beta\right) \right]
 $$
and
 $$
[\psi] = \left[\Ad q_B(W^*S_1W) \circ  \left( \widehat{e_1} \circ \Phi  +_{\beta} \beta\right) \right]
 $$
in $\Ext^{-1}_{J,\beta}(A,B)$. Consider $W$ as a unitary in
$M(IB)$ via the map $\overline{c} : M(B) \to M(IB)$ from the proof
of Lemma \ref{iso42}, and let $\Lambda : A \to Q(IB)$ be the
extension given by
 $$
\Lambda = \Ad q_{IB}(W^*V_1W) \circ \left(\Phi +_{\gamma} \gamma \right) .
 $$
Then $\Lambda$ is equal to $\gamma$ on $J$. By assumption $\Phi$
is invertible which means that it represents an invertible element
of $\Ext(A,IB)$. Hence $\Lambda$ represents also an invertible
element of $\Ext(A,IB)$. As we saw in the proof of Lemma \ref{B35}
this implies that $\Lambda$ represents an invertible element of
$\Ext_{J,\gamma}(A,IB)$. It follows therefore from Lemma
\ref{homtop2} that  $\left[\Ad q_B(W^*S_0W) \circ  \left(
\widehat{e_0} \circ \Phi  +_{\beta} \beta\right) \right] =
{e_0}_*[\Lambda] = {e_1}_*[\Lambda] = \left[\Ad q_B(W^*S_1W) \circ
\left( \widehat{e_1} \circ \Phi  +_{\beta} \beta\right) \right]$
in $\Ext^{-1}_J(A,B)$. Hence $[\psi] =[\varphi]$ in
$\Ext^{-1}_J(A,B)$.
 \end{proof}
\end{prop}

The main virtue of Proposition \ref{iso42}, which we shall exploit
below, is that it gives a description of the relative extension
group without any reference to absorbing $*$-homomorphism.
Furthermore, the description makes it easy to make the relative
extension group functorial, covariantly in the 'coefficient
algebra' $B$, and contravariantly in the pair $J \subseteq A$.

Let $C_{i}$ be the mapping cone of the inclusion $i : J \to A$ which we realize as
 $$
C_{i} = \left\{ f \in IA: \ f(0) = 0, \ f(s) \in J, \ s \in \bigl[\frac{1}{2},1\bigr] \right\} .
 $$
Let $\varphi : A \to Q(B)$ be an invertible extension of $A$ by
$B$ which equals $\alpha$ on $J$. We can then choose a completely
positive contraction $\xi : A \to M(B)$ such that $q_B \circ \xi =
\varphi$. Note that $\xi(j) - \alpha_0(j) \in B$ for any $j\in J$
since $\varphi|_J = q_B \circ \alpha_0$. We define $\Phi: C_i \to
IM(B)$ such that
 $$
\Phi(f)(s) = \begin{cases}  \xi\left(f(s)\right), & s \in \left[0,\frac{1}{2}\right] \\ \left(2-2s\right) \xi\left(f(s)\right) + \left(2s-1\right) \alpha_0(f(s)), & s \in \left[\frac{1}{2},1\right]. \end{cases}
 $$
Then $\Phi$ is a completely positive contraction and $\Phi(fg) -
\Phi(f)\Phi(g) \in SB$ for all $f,g \in C_i$. Thus $\mu(\varphi) =
q_{SB} \circ \Phi: C_i \to Q(SB)$ is an invertible extension of
$C_i$ by $SB$. It is easy to see that we get a group homomorphism
$\mu : \Ext^{-1}_{J,\alpha}(A,B) \to \Ext^{-1}\left(C_i,SB\right)$
such that $\mu[\varphi]= [\mu(\varphi)]$.

\begin{thm}\label{mapiso}
$\mu : \Ext^{-1}_{J,\alpha}(A,B) \to \Ext^{-1}\left(C_i,SB\right)$
is an isomorphism.
\begin{proof} Let $\iota : SA \to C_i$ be the natural embedding, i.e.
 $$
\iota(g)(s) = \begin{cases} g(2s), & s \in
\bigl[0,\frac{1}{2}\bigr]
\\ 0, & s \in \bigl[\frac{1}{2},1\bigr], \end{cases}
 $$
and $p : C_i \to J$ the $*$-homomorphism $p(f) = f(1)$. By
comparing the six-term exact sequence of Theorem \ref{cor1} with
the Puppe sequence of \cite{CS} we see that the five-lemma will
give the theorem if we show that the diagram
\begin{equation}\label{diagiso}
\begin{xymatrix}{
\Ext^{-1}(J,SB) \ar@{=}[d] \ar[r]^-{\partial_0}  &
\Ext^{-1}_{J,\alpha}(A,B) \ar[d]^-{\mu}  \ar[r] &
\Ext^{-1}(A,B)\ar@{=}[d] \\ \Ext^{-1}(J,SB) \ar[r]^-{p^*} &
\Ext^{-1}\left(C_i,SB\right) \ar[r]^-{S^{-1} \circ \iota^*} &
\Ext^{-1}(A,B) }
\end{xymatrix}
\end{equation}
commutes, where $S^{-1}$ is the inverse of the suspension
isomorphism $S : \Ext^{-1}(A,B) \to  \Ext^{-1}(SA,SB)$. To this
end, only the left square requires some care. To prove
commutativity here we consider a unitary $v$ in $M(B)$ such that
$v\alpha_0(j) - \alpha_0(j)v \in B$ for any $j\in J$. Under the
isomorphism $K_1\left(D_{\alpha}(J)\right) \simeq
\Ext^{-1}(J,SB)$, $v$ becomes the extension $\psi = q_{SB} \circ
\psi_0$, where $\psi_0 : J \to IM(B)$ is given by $\psi_0(j)(s) =
(1-s)v \alpha_0(j)v^* + s\alpha_0(j)$. Hence $p^*[v]$ is
represented by the extension $q_{SB} \circ \Phi : C_i \to Q(SB)$,
where $\Phi : C_i \to IM(B)$ is given by
 $$
\Phi(f)(s) = (1-s)v \alpha_0(f(1))v^* + s\alpha_0(f(1)) .
 $$
For comparison $\mu \circ \partial_0[v] $ is represented by the
extension $q_{SB} \circ \Psi : C_i \to Q(SB)$, where $\Psi : C_i
\to  IM(B)$ is given by
 $$
\Psi(f)(s) = \begin{cases}  v\alpha_0\left(f(s)\right)v^*, & s \in \left[0,\frac{1}{2}\right] \\ \left(2-2s\right) v\alpha_0\left(f(s)\right)v^*  + \left(2s-1\right) \alpha_0(f(s)), & s \in \left[\frac{1}{2},1\right]. \end{cases}
 $$
Set $h_{\lambda}(s) = \max \{\lambda, s\}$, and define $\Lambda :
C_i \to  I^2M(B)$ such that
 $$
\Lambda(f)(\lambda,s)  = \begin{cases}  v\alpha_0\left(f\left(h_{\lambda}(s)\right)\right)v^*, & s \in \left[0,\frac{1}{2}\right] \\ \left(2-2s\right) v\alpha_0\left(f\left(h_{\lambda}(s)\right)\right)v^*  + \left(2s-1\right) \alpha_0(\left(h_{\lambda}(s)\right)), & s \in \left[\frac{1}{2},1\right]. \end{cases}
 $$
Then $q_{ISB} \circ \Lambda$ is an extension of $C_i$ by $ISB$
which gives us a homotopy between $q_{SB} \circ \Psi$ and $q_{SB}
\circ \Psi'$, where
 $$
\Psi'(f)(s) = \begin{cases}  v\alpha_0\left(f(1)\right)v^*, & s \in \left[0,\frac{1}{2}\right] \\ \left(2-2s\right) v\alpha_0\left(f(1)\right)v^*  + \left(2s-1\right) \alpha_0(f(1)), & s \in \left[\frac{1}{2},1\right]. \end{cases}
 $$
It is easy to construct a homotopy of invertible extensions
connecting $q_{SB} \circ \Psi'$ to $q_{SB} \circ \Phi$, and we
therefore conclude that the diagram (\ref{diagiso}) commutes.
\end{proof}
\end{thm}

Theorem \ref{mapiso} has many consequences for the relative
extension group. One is that $\Ext^{-1}_{J,\alpha}(A,B) =
\Ext^{-1}(A/J,B)$ when $J$ is a semi-split ideal. Another virtue
of Theorem \ref{mapiso} is that it makes it easy to give the
following description of the map $\partial_0$ from Theorem
\ref{cor1} -- a description which we shall need in Section
\ref{K-hom}.

\begin{lemma}\label{alt}
When the group $\Ext^{-1}(J,SB)$ is identified with $KK(J,B)$ in
the Cuntz picture, and $\Ext^{-1}_{J,\alpha}(A,B)$ is identified
with $\Ext^{-1}_J(A,B)$ via the isomorphism (\ref{iso41}), we have
that $\partial_0\left[\varphi_+|_J,\varphi_-\right] = \left[q_B
\circ \varphi_+,\varphi_-\right]$, where $\varphi_+ : A \to M(B)$
and $\varphi_- : J \to M(B)$ are $*$-homomorphisms such that
$\varphi_+(j) - \varphi_-(j) \in B$ for all $j \in J$.
\begin{proof}
Note that the map (\ref{iso41}) was defined for a particular
absorbing $*$-homomorphism $\beta_0$. Let $u \in M(B)$ be a
unitary such that $\Ad u \circ \alpha_0(a)- \beta_0(a)  \in B$ for
all $a \in A$. There is then an isomorphism
$\Ext_{J,\alpha}^{-1}(A,B) \to \Ext_{J,\beta}^{-1}(A,B)$ defined
by $[\varphi] \mapsto \left[\Ad q_B(u) \circ \varphi\right]$. By
composing with the isomorphism (\ref{iso41}) we obtain an
isomorphism $\nu :\Ext_{J,\alpha}^{-1}(A,B) \to \Ext_J^{-1}(A,B)$.
When $\left(\psi,\psi_0\right)$ is an invertible extension of $A$
by $B$ which splits over $J$ we can define $\Psi : C_i \to IM(B)$
such that
 $$
\Psi(f)(s) = \begin{cases}  \xi\left(f(s)\right), & s \in \left[0,\frac{1}{2}\right] \\ \left(2-2s\right) \xi\left(f(s)\right) + \left(2s-1\right) \varphi_0(f(s)), & s \in \left[\frac{1}{2},1\right], \end{cases}
 $$
where $\xi : A \to M(B)$ is a completely positive contractive lift
of $\varphi$. Then $q_{SB} \circ \Psi$ is an invertible extension
of $C_i$ by $SB$ and we can define a homomorphism $\mu' :
\Ext^{-1}(J,SB)\to\Ext^{-1}(C_i,SB)$ such that
$\mu'\left[\psi,\psi_0\right] = \left[q_{SB} \circ \Psi\right]$.
Then the diagram
\begin{equation}\label{diagiso2}
\begin{xymatrix}{
\Ext^{-1}(J,SB) \ar@{=}[d] \ar[r]^-{\partial_0}  & \Ext^{-1}_{J,\alpha}(A,B) \ar[d]^-{\mu}  \ar[r]^-{\nu} & \Ext^{-1}_J(A,B) \ar[ld]^-{\mu'} \\
\Ext^{-1}(J,SB) \ar[r]^-{p^*} & \Ext^{-1}\left(C_i,SB\right) &  }
\end{xymatrix}
\end{equation}
commutes. The commutativity of the square was established in the
proof of Theorem \ref{mapiso} and it is easy to see that the
triangle commutes. The diagram (\ref{diagiso2}) gives us the lemma
in the following way: The element of $\Ext^{-1}(J,SB)$
corresponding to the Cuntz pair
$\left(\varphi_+|_J,\varphi_-\right)$ is represented by $q_{SB}
\circ \Phi$, where $\Phi : J \to IM(B)$ is given by $\Phi(j)(s) =
(1-s)\varphi_+(j) + s \varphi_-(j)$. Thus $p^*\left[\varphi_+|_J,
\varphi_-\right]$ is represented by $q_{SB} \circ \Lambda$, where
$\Lambda(f)(s) = (1-s)\varphi_+(f(1)) + s \varphi_-(f(1))$. For
comparison $\mu'\left[q_B \circ \varphi_+,\varphi_-\right]$ is
represented by the extension $q_{SB} \circ \Psi'$ where
 $$
\Psi'(f)(s) = \begin{cases}  \varphi_+\left(f(s)\right), & s \in \left[0,\frac{1}{2}\right] \\ \left(2-2s\right) \varphi_+\left(f(s)\right) + \left(2s-1\right) \varphi_-(f(s)), & s \in \left[\frac{1}{2},1\right]. \end{cases}
 $$
Homotopies almost identical to the ones used in the proof of
Theorem \ref{mapiso} now prove that $p^*\left[\varphi_+|_J,
\varphi_-\right] = \mu'\left[q_B \circ
\varphi_+,\varphi_-\right]$. The conclusion of the lemma follows
then from the commutativity in (\ref{diagiso2}) because $\mu'$ is
injective.
\end{proof}
\end{lemma}

Under our standing assumption that there is an absorbing
$*$-homorphism $A \to M(B)$ whose restriction to $J$ is also
absorbing, every element of $KK(J,B)$ can be represented by a
Cuntz pair of the form considered in Lemma \ref{alt}.

\section{When $A$ and $J$ have the same unit.}\label{unit}

Assume now that the pair $J \subseteq A$ share the same unit $1
\in J$. It is then often natural and convenient to consider
extensions that are unital. This is certainly the case for the
relative K-homology of compact spaces which we investigate in the
following section. In the present section we describe how to
adjust the definitions and results of the previous sections in
order to 'fix the unit'. Most of the considerations are standard,
so we will be brief (but, hopefully, precise).

First of all the role of the absorbing $*$-homomorphisms must now
be taken by the \emph{unitally} absorbing $*$-homomorphisms. The
first lemma shows that this does not effect the fundamental
condition of Section \ref{sixterms}.

\begin{lemma}\label{unit1}
There is an absorbing $*$-homomorphism $\alpha_0 : A \to M(B)$
such that ${\alpha_0}|_J: J \to M(B)$ is also absorbing if and
only if there is a unitally absorbing $*$-homomorphism $\beta_0 :
A \to M(B)$ such that $\beta_0|_J: J \to M(B)$ is also unitally
absorbing.
\begin{proof}
Assume first that $\alpha_0$ exists. Since there exists a unitally
absorbing $*$-homomorphism $A \to M(B)$, \cite{Th1}, it follows
from Lemma 1.1 in \cite{MT} that there is an absorbing
$*$-homomorphism $A \to M(B)$ such that the image of $1$ is the
range projection of an isometry in $M(B)$, and then it follows
from the essential uniqueness of absorbing $*$-homomorphisms that
this is the case for all of them. In particular, there is an
isometry $W \in M(B)$ such that $WW^* = \alpha_0(1)$. Then
$W^*\alpha_0(\cdot)W : A \to M(B)$ is a unital $*$-homomorphism
and we claim that it is unitally absorbing. To show this we check
that condition 1) of Theorem 2.1 of \cite{Th1} is fullfilled. Let
$\varphi : A \to B$ be a completely positive contraction. Extend
$\varphi$ to $A^+ = A \oplus \mathbb C$ by annihilating the direct
summand $\mathbb C$. Since $\alpha_0^+ : A^+ \to M(B)$ is unitally
absorbing by Theorem 2.7 of \cite{Th1}, there is a sequence
$\{W_n\} \subseteq M(B)$ such that $\lim_{n \to \infty} W_n^*b =
0$ for all $b \in B$ and $\lim_{n \to \infty} W_n^*\alpha_0(a)W_n
= \varphi(a)$ for all $a \in A$. It follows that $\lim_{n \to
\infty} W_n^*Wb = 0$ for all $b \in B$ and $\lim_{n \to \infty}
W_n^*WW^*\alpha_0(a)WW^*W_n = \varphi(a)$ for all $a$, verifying
that $W^*\alpha_0(\cdot)W$ is unitally absorbing. The same
argument applies to its restriction to $J$, and hence $\beta_0 =
W^*\alpha_0(\cdot)W$ is unitally absorbing on both $A$ and $J$.

Conversely, when $\beta_0 : A \to M(B)$ is unitally absorbing on
both $A$ and $J$, Lemma 1.1 of \cite{MT} shows that there is an
isometry $V \in M(B)$ such that $\alpha_0 = \Ad V \circ \beta_0$
is absorbing on both $A$ and $J$.
\end{proof}
\end{lemma}

Assume now that $\beta_0 : A \to M(B)$ is a unitally absorbing
$*$-homomorphism such that $\beta_0|_J : J \to M(B)$ is also
unitally absorbing. Set $\beta = q_B \circ \beta_0$. We say that
two unital extensions, $\varphi, \psi : A \to Q(B)$, that are
equal to $\beta$ on $J$, are unitarily equivalent when there is a
unitary connected to 1 in the unitary group of $\beta(J)'\cap
Q(B)$ such that $\Ad v \circ \varphi = \psi$. By repeating the
arguments that proved Lemma \ref{B94} we find that the addition
$+_{\beta}$, defined using isometries from $\beta(A)' \cap Q(B)$,
gives the unitary equivalence classes of unital extensions of $A$
by $B$ which equal $\beta$ on $J$ the structure of an abelian
semi-group. We define stable equivalence in this setting in the
natural way: $\varphi$ and $\psi$ are stably equivalent, as unital
extensions which equal $\beta$ on $J$, when $\varphi +_{\beta}
\beta$ is unitarily equivalent to $\psi +_{\beta} \beta$. The
stable equivalence classes of unital extensions of $A$ by $B$
which equal $\beta$ on $J$ is then an abelian semi-group with $0$,
and the invertible elements of this semi-group form an abelian
group which we denote by $\Ext_{J,\beta}^{-1}(A,B)$. Let $V, W$ be
isometries in $M(B)$ such that $VV^* + WW^* = 1$. By Lemma 1.1 of
\cite{MT}, $\alpha_0 = \Ad V \circ \beta_0$ will be absorbing on
both $A$ and $J$ and we can define a group homomorphism
$\chi_0:\Ext_{J,\beta}^{-1}(A,B) \to \Ext^{-1}_{J,\alpha}(A,B)$
such that $\chi_0[\varphi] = [\Ad q_B(V) \circ \varphi]$. In the
other direction, if $\psi : A \to Q(B)$ is an extension which
equals $\alpha$ on $J$, note that $\Ad q_B(V)^* \circ\psi$ is a
unital extension of $A$ by $B$ which equals $\beta$ on $J$. We
define a homomorphism $\chi_1 : \Ext^{-1}_{J,\alpha}(A,B)   \to
\Ext_{J,\beta}^{-1}(A,B)$ such that $\chi_1[\psi] = [\Ad q_B(V)^*
\circ \psi]$. It is easy to see that $\chi_1$ is the inverse of
$\chi_0$. Hence we see that the unital version of the relative
extension group agrees with non-unital version.

Note that in the particular case where $A$ is unital and $J =
\mathbb C1 \subseteq A$, the group $\Ext_{J,\beta}^{-1}(A,B)$ is
the group $\Ext^{-1}_{unital}(A,B)$ considered in \cite{MT}, and
the six-term exact sequence of Theorem \ref{cor1} is then the
six-term exact sequence of Skandalis, \cite{S}, a construction of
which was exhibited in \cite{MT}. In fact, in the present setting
where $J$ and $A$ have a common unit the six-term exact sequence
of Theorem \ref{cor1} can be modified so that the involved
extension groups are 'unital' in the sense that they are defined
using unital extensions only. The key point for this is that since
$\beta_0$ is unitally absorbing, the isomorphism (\ref{B51}) can
be substituted by an isomorphism $K_1\left(\beta(A)^{\prime} \cap
Q(B)\right) \simeq KK(A,B)$ because $\mathcal
D_{\alpha}(A)/\mathcal X_{\alpha}(A) \simeq \beta(A)^{\prime} \cap
Q(B)$, cf. \cite{MT}. By using this isomorphism in place of
(\ref{B51}) and the isomorphism $K_1\left(\beta(J)^{\prime} \cap
Q(B)\right) \simeq KK(J,B)$ in place of (\ref{B32}), the proof of
Theorem \ref{cor1} can easily be adopted to yield the following
six-term exact sequence in the present case:
\begin{equation*}\label{nycor1}
\begin{xymatrix}{
 \Ext^{-1}_J(A,B)  \ar[r]  &   \Ext^{-1}_{unital}(A,B)  \ar[r]^-{i^*} & \Ext^{-1}_{unital}(J,B)  \ar[d]^-{\partial_1} \\
\Ext^{-1}_{unital}(J,SB)  \ar[u]^-{\partial_0} &   \Ext^{-1}_{unital}\left(A,SB\right) \ar[l]^-{i^*} &  \Ext^{-1}_J(A,SB) \ar[l] }
\end{xymatrix}
\end{equation*}

As one would expect by now, the alternative picture of the
relative extension group given in Section \ref{other} is also not
changed by restricting attention entirely to unital extensions.
This will be very useful in the following, so let us make it
precise: By a unital extension of $A$ by $B$ which splits over $J$
we mean a pair $(\varphi, \varphi_0)$ where $\varphi : A \to Q(B)$
is a unital extension of $A$ by $B$, and $\varphi_0 : J \to M(B)$
is a unital $*$-homomorphism such that $q_B \circ \varphi_0 =
\varphi|_J$. Recall that if $\varphi$ is invertible, it is
actually unitally invertible, i.e. there is another unital
extension $\psi$ of $A$ by $B$ with the property that $\varphi
\oplus \psi$ is a split extension (of $A$ by $B$). Two unital
invertible extensions, $(\varphi,\varphi_0)$ and $(\psi,\psi_0)$,
of $A$ by $B$ which split over $J$ are now homotopic in norm when
there is a path $\left(\Phi^t,\Phi_0^t\right), t \in [0,1]$, of
unital invertible extensions of $A$ by $B$ which split over $J$
such that $[0,1] \ni t \mapsto \Phi^t(a)$ is norm-continuous for
all $a \in A$ and $[0,1] \ni t \mapsto \Phi_0^t(j)$ is
norm-continuous for all $j \in J$, $\left(\Phi^0,\Phi_0^0\right) =
(\varphi,\varphi_0)$ and $\left(\Phi^1,\Phi_0^1\right) =
(\psi,\psi_0)$. We say that $(\varphi,\varphi_0)$ and
$(\psi,\psi_0)$ are stably homotopic in norm when there is a
unital $*$-homomorphism $\pi : A \to M(B)$ such that $\left(
\varphi \oplus q_B \circ \pi, \varphi_0 \oplus \pi|_J\right)$ and
$\left( \psi \oplus q_B \circ \pi, \psi_0 \oplus \pi|_J\right)$
are homotopic in norm. The stable homotopy classes of unital
invertible extensions of $A$ by $B$ which split over $J$ form a
group which we temporarily denote by $\Ext_J^{-1,+}(A,B)$.

\begin{lemma}\label{forgetful}
The forgetful map $\Ext^{-1,+}_J(A,B) \to \Ext^{-1}_J(A,B)$ is an
isomorphism.
\begin{proof}
Surjectivity: Let $\left(\varphi,\varphi_0\right)$ be an
invertible extension of $A$ by $B$ which splits over $J$. By
adding $(0,0)$ we don't change the class of
$\left(\varphi,\varphi_0\right)$, but reach a situation where
there is a $*$-homomorphism $\pi : A \to M(B)$ such that $\pi(1) =
\varphi_0(1)^{\perp}$. Then, with $\overline{\Theta}$ the
$*$-isomorphism from the proof of Lemma \ref{B73}, $w =
\overline{\Theta}\left( \begin{smallmatrix} \varphi_0(1)  & \pi(1)
\\ 0 & 0 \end{smallmatrix} \right)$ is a partial isometry such
that $\Ad q_B(w) \circ \left( \varphi \oplus q_B \circ \pi\right)
$ and $\Ad w \circ \left( \varphi_0 \oplus \pi|_J\right)$ are both
unital. By choosing a unitary dilation of $w$ and using that the
unitary group of $M(B)$ is connected in the norm-topology,
\cite{M}, \cite{CH}, we see that the class of
$\left(\varphi,\varphi_0\right)$ in $\Ext^{-1}_J(A,B)$ is also
represented by the unital pair $\left( \Ad q_B(w) \circ \left(
\varphi \oplus q_B \circ \pi\right),\Ad w \circ \left( \varphi_0
\oplus \pi|_J\right) \right)$.

Injectivity: Let $\left(\varphi,\varphi_0\right)$ and
$\left(\psi,\psi_0\right)$ be unital and invertible extensions of
$A$ by $B$ which split over $J$ and define the same element of
$\Ext^{-1}_J(A,B)$. After the addition of a pair $\left(q_B \circ
\pi, \pi|_J\right)$ we have a homotopy in norm, $\left(\Phi_t,
\Phi^0_t\right), t \in [0,1]$, connecting  $\left(\varphi \oplus
q_B \circ \pi,\varphi_0 \oplus \pi|_J \right)$ to $\left(\psi
\oplus q_B \circ \pi,\psi_0 \oplus \pi|_J\right)$. Standard
arguments show that there is a normcontinuous path $U_t, t \in
[0,1]$, of unitaries in $M(B)$ such that $\Phi^0_t(1) =U_tpU_t^*$
for all $t \in [0,1]$, where $p = 1 \oplus \pi(1)$. Note that
$U_0pU_0^* = U_1pU_1^* = p$. By the same trick of adding $(0,0)$
as above, we can arrange that there is a $*$-homomorphism $\chi: A
\to M(B)$ such $\chi(1) = p^{\perp}$. Define $\Psi^0_t(j) =
\Phi^0_t(j) + \Ad U_t \circ \chi(j)$ and $\Psi_t(a) = \Phi_t(a)  +
\Ad q_B(U_t) \circ q_B \circ \chi(a)$. Then $\left(\Psi_t,
\Psi^0_t\right)$ is a homotopy in norm connecting
 $$
\left( \left(\varphi \oplus q_B\circ \pi\right) + q \circ \Ad U_0
\circ \chi, \left(\varphi_0 \oplus \pi|_J\right) + \Ad U_0 \circ
\chi|_J \right)
 $$
to
 $$
\left( \left(\psi \oplus q_B\circ \pi\right) + q_B \circ \Ad U_1
\circ \chi, \left(\psi_0 \oplus \pi|_J\right) + \Ad U_1 \circ
\chi|_J \right).
 $$
Since there are unital $*$-homomorphisms $\pi_{\pm} : A \to M(B)$
such that
\begin{equation*}
\begin{split}
&\left( \left(\varphi \oplus q_B\circ \pi\right) + q_B \circ \Ad U_0 \circ \chi, \left(\varphi_0 \oplus \pi\right) + \Ad U_0 \circ \chi|_J \right) \\
&= \left( \varphi \oplus q_B \circ \pi_+, \varphi_0 \oplus \pi_+|_J\right)
\end{split}
\end{equation*}
and
\begin{equation*}
\begin{split}
&\left( \left(\psi \oplus q_B\circ \pi\right) + q_B \circ \Ad U_1 \circ \chi, \left(\psi_0 \oplus \pi\right) + \Ad U_1 \circ \chi|_J \right) \\
&= \left( \psi \oplus q_B \circ \pi_-, \psi_0 \oplus \pi_-|_J\right),
\end{split}
\end{equation*}
we conclude that $\left( \varphi, \varphi_0\right)$ and $\left(
\psi, \psi_0\right)$ define the same element of
$\Ext^{-1,+}_J(A,B)$.
\end{proof}
\end{lemma}

Lemma \ref{forgetful} is our excuse for not distinguishing between
$\Ext^{-1,+}_J(A,B)$ and $\Ext^{-1}_J(A,B)$ in the following.

\section{Relative K-homology for spaces.}\label{K-hom}

Fix a separable infinite-dimensional Hilbert space $H$ and denote
by $\mathbb L(H)$ the algebra of bounded operators on $H$, and by
$\mathbb K$ the ideal in $\mathbb L(H)$ consisting of the compact
operators. In this section we will study the relative extension
group $\Ext_{J}^{-1}\left(A, B\right)$ in the case when the
'coefficient algebra' is $B = \mathbb K$ and $J \subseteq A$ is a
unital inclusion of one abelian $C^*$-algebra into another. Since
the coefficients are now fixed we drop them in the notation. In
the same vein we write $Q$ for the Calkin algebra and $q : \mathbb
L(H) \to Q$ for the quotient map. In order to draw directly on the
work of Brown, Douglas and Fillmore, \cite{BDF2}, we will use the
'unital version' of the relative extension group, as explained in
Section \ref{unit}.

\begin{lemma}\label{BDF-1}
Let $A \subseteq \mathbb L(H)$ be a separable $C^*$-subalgebra
such that $q\left(A\right) \subseteq Q$ is abelian. Let $\omega_i
: q\left(A\right) \to \mathbb C, i = 1,2, \dots $, be a sequence
of characters of $q\left(A\right)$, and let $\left\{a_i\right\}$
be a dense sequence in $A$.

There is then a family of continuous maps $\chi_i : [1,\infty) \to
H, i = 1,2, 3, \dots,$ such that
\begin{enumerate}
\item[i)] $\left( \chi_1(t),\chi_2(t),\chi_3(t), \dots \right)$ is an orthonormal set for all $t \in [1,\infty)$,
\item[ii)] $\chi_i(t) = \chi_i(s), i \geq n$, when $s,t \in [1,n]$,
\item[iii)] $\sum_{i=1}^{\infty}\left\| a_i \chi_i(t) -  \omega_i\left(q(a_i)\right) \chi_i(t)\right\|^2 < \infty$ for all $t \in [1,\infty)$,
\item[iv)] $\lim_{ t \to \infty} \sup_i \left\| a\chi_i(t) - \omega_i\left(q(a)\right) \chi_i(t)\right\| = 0$ for all $a \in A$, and
\item[v)] $\lim_{t \to \infty} \sup_i \left\| k\chi_i(t)\right\| = 0$ for all $k \in \mathbb K$.
\end{enumerate}
\begin{proof}
Let $\{\mu_i\}$ be a sequence of characters on $q\left(A\right)$
with the property that each $\omega_j$ occurs infinitely many
times in $\{\mu_i\}$, and let $b_1,b_2, b_3, \dots$, be a dense
sequence in $\mathbb K$. By Lemma 1.3 of \cite{BDF2} there is an
orthonormal sequence $\psi_i, i = 1,2, \dots$, in $H$ such that
 $$
\left\|a_k \psi_i - \mu_i\left(q(a_k)\right)\psi_i\right\| \leq 2^{-i}
 $$
when $k \leq i$. Since each $\omega_j$ occurs infinitely often in
$\{\mu_i\}$ we can select subsequences
$\left\{\varphi_i^k\right\}_{i=1}^{\infty}, k = 1,2,3, \dots$,
from $\left\{\psi_i\right\}$ such that
\begin{enumerate}
\item[a)] $\left\{\varphi_i^k\right\}_{i=1}^{\infty} \cap \left\{\varphi_i^{k'}\right\}_{i=1}^{\infty} = \emptyset$ when $k \neq k'$, and
\item[b)] $\sup_i \left\| a_j \varphi_i^k - \omega_i\left(q\left(a_j\right)\right)\varphi_i^k \right\| \leq \frac{1}{k}$,
\item[c)] $\sup_i \left\|b_j \varphi_i^k \right\| \leq \frac{1}{k}$,
\end{enumerate}
for all $j =1,2, \dots, k$. Set
 $$
\gamma_i^k = \begin{cases} \varphi_i^k, & i \leq k \\ \varphi_i^1, & i > k .\end{cases}
 $$
and
 $$
\chi_i(t) = \sqrt{t -n}\gamma_i^{n+1} + \sqrt{n+1 -t} \gamma_i^n
 $$
when $t \in [n,n+1]$. It is straigthforward to show that $\{\chi_i\}$ has the properties i)-v).
\end{proof}
\end{lemma}

Let $X$ and $Y$ be compact Hausdorff spaces and $f : X \to Y$ a
continuous surjection. There is then a unital embedding $i : C(Y)
\to C(X)$ such that $i(g) = g \circ f$, and we will sometimes
identify $C(Y)$ with $i\left(C(Y)\right) \subseteq C(X)$. It
follows from \cite{K2} that any unital $*$-homomorphism $\alpha_0
: C(X) \to \mathbb L(H)$ such that $\alpha = q \circ \alpha_0$ is
injective, is also unitally absorbing. In particular, it follows
that $\alpha_0|_{C(Y)}$ is unitally absorbing whenever $\alpha_0$
is. We fix a unitally absorbing $*$-homomorphism $\alpha_0 : C(X)
\to \mathbb L(H)$, and denote $\Ext_{C(Y),\alpha}(C(X),\mathbb K)$
by
 $$
\Ext_{Y,f}(X).
 $$
It's this group we shall investigate in this section.

\begin{lemma}\label{BDF11}
Let $\varphi : C(X) \to Q$ be an injective and unital extension of
$C(X)$ by $\mathbb K$. Let $\varphi' : C(X) \to \mathbb L(H)$ be a
continuous map such that $q_{\mathbb K} \circ \varphi^{\prime} =
\varphi$. There is then a continuous path $V_t, t \in [1,
\infty)$, of unitaries in $\mathbb L(H)$ such that
\begin{enumerate}
\item[i)] $\Ad q\left(V_t\right) \circ  \left( \varphi \oplus \alpha\right) = \varphi$ for all $t \in [1,\infty)$,
\item[ii)] $\lim_{t \to \infty} V_t \left( \varphi' \oplus \alpha_0\right)(g)V_t^* = \varphi'(g)$ for all $g \in C(X)$,
\item[iii)] $V_1 - V_t \in \mathbb K$ for all $t \in [1,\infty)$.
\end{enumerate}
\begin{proof}
We pick a dense sequence $\{x_i\}_{i=1}^{\infty}$ in $X$ and an
orthonormal basis $\{e_i\}_{i=1}^{\infty}$ in $H$. We arrange that
each point $x_j$ is repeated infinitely many times in
$\{x_i\}_{i=1}^{\infty}$. Define $D : C(X) \to \mathbb L(H)$ such
that $D(g)\psi = \sum_{i=1}^{\infty} g(x_i)\left<\psi,
e_i\right>e_i$. Apply then Lemma \ref{BDF-1} with $A$ the
$C^*$-algebra generated by $\varphi^{\prime}\left(C(X)\right)$ and
$\omega_i(a) =  \varphi^{-1}(a)(x_i)$ to get the continuous
functions $\chi_i : [1,\infty) \to H$ with the properties listed
there. Let $\sigma_1, \sigma_2 : \mathbb N \to \mathbb N$ be
injective maps such that $\sigma_1\left(\mathbb N\right) \cap
\sigma_2\left(\mathbb N\right) = \emptyset$, $\mathbb N =
\sigma_1\left(\mathbb N\right) \cup \sigma_2\left(\mathbb
N\right)$ and $x_i = x_{\sigma_1(i)} =  x_{\sigma_2(i)}$ for all
$i \in \mathbb N$ (such maps exist because each $x_j$ is repeated
infinitely many times).

Set $P_t\psi = \sum_{i=1}^{\infty} \left< \psi, \chi_i(t)\right> \chi_i(t)$. Then $P_t$ is a projection and we define $V_t: P_tH \to H$ such that $V_t \chi_i(t) = \chi_{\sigma_1(i)}(t)$. Define isometries $S_t,T_t : H \to H$ such that $S_t = P_t^{\perp} + V_tP_t$ and $T_te_i = \chi_{\sigma_2(i)}(t)$. It follows from the properties i)-v) of $\{\chi_i\}$ that both $t \mapsto S_t$ and $t \mapsto T_t$ are continuous in norm, that $S_t\varphi^{\prime}(g)S_t^* + T_t D(g)T_t^* - \varphi^{\prime}(g) \in \mathbb K$ and that $\lim_{t \to \infty}S_t\varphi^{\prime}(g)S_t^* + T_t D(g)T_t^* = \varphi^{\prime}(g)$. Let $V_1,V_2$ be the isometries used to define the addition $\oplus$, and set $U_t = S_tV_1^* + T_tV_2^*$. Then
\begin{enumerate}
\item[a)] $U_t \left( \varphi^{\prime} \oplus D\right)(g)U_t^* - \varphi^{\prime}(g) \in \mathbb K$ for all $t \in [1,\infty)$ and all $g \in C(X)$,
\item[b)] $\lim_{t \to \infty} U_t \left( \varphi^{\prime} \oplus D\right)(g)U_t^* = \varphi^{\prime}(g)$ for all $g \in C(X)$, and
\item[c)] $U_1 - U_t \in \mathbb K$ for all $t \in [1, \infty)$.
\end{enumerate}
By Theorem 3.11 of \cite{DE2} there is is a norm-continuous path
$W_t, t \in [1,\infty)$, of unitaries in $\mathbb L(H)$ such that
$\Ad W_t \circ D(g) - \alpha_0(g) \in \mathbb K$ for all $t \in
[1,\infty), g \in C(X)$, and $\lim_{t \to \infty} \Ad W_t \circ
D(g) = \alpha_0(g)$ for all $g \in C(X)$. It follows then from
Lemm \ref{iso6} that the Cuntz-pair $\left(\Ad W_1 \circ D,
\alpha_0\right)$ represents zero in $KK(C(X),\mathbb K)$, and then
an application of Theorem 3.12 of \cite{DE2} shows that we can
assume that $W_t - W_1 \in \mathbb K$ for all $t$. Set $V_t =
U_t\left(1 \oplus W_t^*\right)$.
\end{proof}
\end{lemma}

\begin{thm}\label{triv}
Let $\varphi : C(X) \to Q(\mathbb K)$ be an injective and unital
extension of $C(X)$ by $\mathbb K$. Assume that there is a unital
$*$-homomorphism $\varphi_0 : C(Y) \to \mathbb L(H)$ such that
$\varphi \circ i = q \circ \varphi_0$. Then the following are
equivalent:
\begin{enumerate}
\item[i)] $(\varphi,\varphi_0)$ represents zero in $\Ext_{Y,f}(X)$.

\item[ii)] There is a path $\psi_t, t \in [1,\infty)$, of
$*$-homomorphisms $\psi_t : C(X) \to \mathbb L(H)$ such that $t
\mapsto \psi_t(g)$ is continuous for all $g \in C(X)$, $\varphi =
q \circ \psi_t$ for all $t \in [1,\infty)$, and $\lim_{t \to
\infty} \psi_t \circ i(h) = \varphi_0(h)$ for all $h \in C(Y)$.

\item[iii)] There is a sequence of $*$-homomorphisms
$\psi_n : C(X) \to \mathbb L(H)$ such that $q \circ \psi_n =
\varphi$ for all $n$, and $\lim_{n \to \infty} \psi_n \circ i(h) =
\varphi_0(h)$ for all $h \in C(Y)$.
\end{enumerate}
\begin{proof}
The implication ii) $\Rightarrow$ iii) is trivial so it suffices
to prove that i) $\Rightarrow$ ii) and that iii) $\Rightarrow$ i).
First i) $\Rightarrow$ ii): It follows from Section \ref{unit}
that there is a unital $*$-homomorphism $\pi : C(X) \to \mathbb
L(H)$ such that $\left( \varphi \oplus q_B \circ \pi, \varphi_0
\oplus \pi|_{C(Y)}\right)$ is homotopic in norm, as a unital
extension of $C(X)$ by $\mathbb K$ which splits over $C(Y)$, to
the pair $\left( q_B \circ \pi, \pi|_{C(Y)}\right)$. It follows
then from \cite{BDF2} that $\varphi$ represents zero in $\Ext(X)$.
There is therefore a unital $*$-homomorphism $\varphi' : C(X) \to
\mathbb L(H)$ such that $q \circ \varphi' = \varphi$. Then
$(\varphi' \circ i, \varphi_0)$ is a Cuntz pair and from the
description of $\partial$ given after Proposition \ref{iso42}, we
find that  $\partial \left[\varphi' \circ i,\varphi_0\right]
=\left[\varphi,\varphi_0\right]$. Since
$\left[\varphi,\varphi_0\right] = 0$, the six-term exact sequence
of Theorem \ref{cor1} shows that $ \left[\varphi' \circ
i,\varphi_0\right] = i^*\left[ \psi_+, \psi_-\right]$, where
$\psi_{\pm}: C(X) \to \mathbb L(H)$ are (not necessarily unital)
$*$-homomorphisms such that $\psi_+(g) - \psi_-(g) \in \mathbb K$
for all $g \in C(X)$. By adding the same $*$-homomorphism to
$\psi_+$ and $\psi_-$ we may assume that $q \circ \psi_{\pm}$ are
both injective. Since $\left[\varphi' \circ i \oplus \psi_- \circ
i, \varphi_0 \oplus \psi_+ \circ i\right] = 0$ in $KK(C(Y),\mathbb
K)$ and $\varphi'$ and $\varphi_0$ are both unital it follows that
$\left[\psi_+(1),\psi_-(1)\right] = 0$ in $KK(\mathbb C,\mathbb
K)$. Thus also $\left[1- \psi_+(1),1-\psi_-(1)\right] = 0$ in
$KK(\mathbb C,\mathbb K)$. Let $\chi : C(X) \to \mathbb C$ be any
character. Then $\left[\chi \left(1 - \psi_+(1)\right), \chi
\left(1 - \psi_-(1)\right)\right] = 0$ in $KK\left(C(X), \mathbb
K\right)$. It follows that $\left[ \psi_+, \psi_-\right] = \left[
\psi_+, \psi_-\right] + \left[\chi \left(1 - \psi_+(1)\right),
\chi \left(1 - \psi_-(1)\right)\right]$ can be represented by a
Cuntz pair $\gamma_{\pm}$ of unital $*$-homomorphisms
$\gamma_{\pm} : C(X) \to \mathbb L(H)$ such that $q \circ
\gamma_+$ and $q\circ \gamma_-$ are both injective. Since
$\left[\varphi' \circ i \oplus \gamma_- \circ i, \varphi_0 \oplus
\gamma_+ \circ i\right] = 0$ in $KK(C(Y),\mathbb K)$, it follows
from Theorem 3.12 of \cite{DE2} there is a continuous path $V_t, t
\in [1,\infty)$, of unitaries in $1 + \mathbb K$ such that
 $$
\lim_{t \to \infty} \Ad V_t  \left(\varphi' \circ i \oplus \gamma_- \circ i\right) = \varphi_0 \oplus \gamma_+ \circ i .
 $$
Set $\gamma = q \circ \gamma_+ = q\circ \gamma_-$. It follows from
Lemma \ref{BDF11} that there are paths of unitaries $U_t, W_t$,
$t\in[1,\infty)$, such that
 $$
\Ad q\left(W_t\right) \left( \varphi \oplus \gamma\right) = \varphi = \Ad q\left(U_t\right) \left( \varphi \oplus \gamma\right)
 $$
for all $t$,
 $$
\lim_{t \to \infty} \Ad W_t \left( \varphi' \oplus \gamma_-\right)(g) = \varphi'(g)
 $$
for all $g \in C(X)$ and
 $$
\lim_{t \to \infty} \Ad U_t \left( \varphi_0 \oplus \gamma_+ \circ i\right)(h) = \varphi_0(h)
 $$
for all $h \in C(Y)$. Set $T_t = U_tV_tW_t^*$ and $\psi_t = \Ad T_t \circ \varphi'$.

iii) $\Rightarrow$ i): Let $\Psi : C(X) \to \mathbb L(H)$ be the
$*$-homomorphism obtained from the representation $\diag \left(
\psi_1, \psi_2, \psi_3, \dots \right)$. Then $\left( q \circ \Psi,
\Psi \right)$ represents zero in  $\Ext_{Y,f}\left(X\right)$. Let
$\varphi' : C(X) \to \mathbb L(H)$ be a continuous function such
that $q \circ \varphi' = \varphi$ and $\varphi' \circ i =
\varphi_0$. Let $\chi_t, t \in [1,\infty)$, be the path of maps
such that $\chi_t, t \in [n,n+1]$, connects
 $$
\diag( \psi_1, \psi_2, \dots, \psi_{n-1},\varphi',\psi_n,\psi_{n+1}, \dots \dots )
 $$
to
 $$
\diag( \psi_1, \psi_2, \dots, \psi_{n},\varphi',\psi_{n+1}, \psi_{n+2}, \dots \dots )
 $$
by rotation. Considered as maps $\chi_t : C(X) \to \mathbb L(H)$
they give us a path of maps such that $q \circ \chi_t = \varphi
\oplus  q \circ \Psi$ while $\chi_t \circ i$ is a $*$-homomophism,
and $\lim_{t \to \infty} \chi_t \circ i(h) = \Psi(h)$ for all $h
\in C(Y)$. It follows that $\left( \varphi, \varphi_0\right)
\oplus \left( q \circ \Psi, \Psi \right)$ is homotopic in norm to
$\left( q \circ \Psi, \Psi \right)$.
\end{proof}
\end{thm}

Two injective unital extensions $\varphi, \psi : C(X) \to Q$ that
equal $\alpha$ on $C(Y)$ are said to \emph{equivalent} when there
is a norm-continuous path $V_t, t \in [1, \infty)$, of unitaries
in $\mathbb L(H)$ such that
\begin{enumerate}
\item[i)] $\Ad q\left(V_t\right) \circ \varphi = \psi$ for all $t \in [1, \infty)$,
\item[ii)] $\lim_{t \to \infty} V_t\alpha_0 \circ i(h)V_t^* = \alpha_0 \circ i(h)$ for all $h \in C(Y)$, and
\item[iii)] $V_1 - V_t \in \mathbb K$ for all $t \in [1, \infty)$.
\end{enumerate}
We write $\varphi \simeq \psi$ in this case.

\begin{thm}\label{BDF12}
Let $\varphi, \psi : C(X) \to Q$ be injective unital extensions
that equal $\alpha$ on $C(Y)$. Then the following are equivalent:
\begin{enumerate}
\item[i)] $\varphi$ and $\psi$ define the same element of $\Ext_{Y,f}\left(X\right)$.
\item[ii)] $\varphi \ \simeq \ \psi$.
\item[iii)] There is a unitary $V$ connected to $1$ in the unitary group of a
 $$
D=\left\{ m \in \mathbb L(H): \ m\alpha_0 \circ i(h) -
\alpha_0\circ i(h)m \in \mathbb K , \ h \in C(Y) \right\}
 $$
such that $\Ad q(V) \circ \varphi = \psi$.

\end{enumerate}

\begin{proof} It is trivial that iii) implies i).

i) $\Rightarrow$ ii): Assume that $[\varphi] = [\psi]$ in
$\Ext_{Y,f}\left(X\right)$. There is then a unitary $V$ connected
to $1$ in the unitary group of $D$ such that
 $$
\Ad q (V) \circ \left( \varphi +_{\alpha} \alpha\right) =  \psi +_{\alpha} \alpha .
 $$
In particular, it follows that $\left(\beta_1, \beta_2\right) =
\left(\Ad V \circ \alpha_0 \circ i, \alpha_0 \circ i\right)$ is a
Cuntz pair, and since $V$ is connected to $1$ in the unitary group
of $D$, the pair represents zero in $KK(C(Y),\mathbb K)$. It
follows therefore from Theorem 3.12 of \cite{DE2} that there is a
norm-continuous path $S_t, t \in [1,\infty)$, of unitaries in
$1+\mathbb K$ such that $\lim_{t \to \infty} S_t\beta_1(h)S_t^* -
\beta_2(h) = 0$ for all $h \in C(Y)$. It follows that $U_t = S_tV,
t \in [1,\infty)$, is a norm-continuous path of unitaries in
$\mathbb L(H)$ giving rise to an equivalence $\varphi +_{\alpha}
\alpha \ \simeq \ \psi +_{\alpha} \alpha$. By Lemma \ref{BDF11}
$\varphi +_{\alpha} \alpha$ and $\psi +_{\alpha} \alpha$ are
equivalent to $\varphi$ and $\psi$, respectively. Thus $\varphi$
and $\psi$ are equivalent.

ii) $\Rightarrow$ iii): Let $V_t, t\in [1,\infty)$, be a
continuous path of unitaries in $\mathbb L(H)$ giving rise to the
equivalence $\varphi \ \simeq \psi$. It suffices to show that
$V_1$ is connected to $1$ in the unitary group of $D$. Note that
the Cuntz pair $\left(\Ad V_1 \circ \alpha_0, \alpha_0\right)$
represents $0$ in $KK(C(Y), \mathbb K)$ by Lemma \ref{iso6}. It
follows therefore from Paschke's duality result, \cite{Pa}, that
$q(V_1)$ is connected to $1$ in the unitary group of the relative
commutant $\left(q\circ \alpha_0 \circ i(C(Y))\right)^{\prime}
\cap Q$. It follows that there is a unitary $W$ connected to $1$
in the unitary group of $D$ such that $q\left(W\right) =
q\left(V_1\right)$. Then $V_1W^* \in 1 + \mathbb K$ and since the
unitary group of $1 + \mathbb K$ is connected, we see that also
$V_1$ is connected to $1$ in the unitary group of $D$.

\end{proof}
\end{thm}

It follows from Theorem \ref{BDF12} that $\Ext_{Y,f}(X)$ is the
group of equivalence classes of injective (or essential) unital
extensions of $C(X)$ by $\mathbb K$ that equal $\alpha$ on $C(Y)$
with the addition defined by (\ref{B1}).

\section{Normal operators.}
In this section we use the results of the previous sections to
prove Theorem \ref{introthm} from the introduction. The key point
is the following

\begin{thm}\label{ann6}
Let $X$ be a compact metric space and $Y$ a compact subset of the
complex plane $\mathbb C$. Let $f : X \to Y$ be a continuous
surjection. Then the map $f_* : K_0(X) \to K_0(Y)$ is surjective.
\end{thm}

To prove this recall that for every compact metric space $Y$ there
is a natural decomposition $K_0\left(C(Y)\right) = \tilde{K}^0(Y)
\oplus C(Y,\mathbb Z)$, where the summand $\tilde{K}^0(Y)$ is
called the reduced K-theory of $Y$, at least when $Y$ is
connected. We say that $Y$ \emph{has trivial reduced K-theory}
when $\tilde{K}^0(Y) =0$. It is well-known and easy to see that a
compact subset of the complex plane has trivial reduced K-theory
and trivial $K^1$-group. Therefore Theorem \ref{ann6} is a
consequence of the following more general result.

\begin{thm}\label{ann7}
Let $X$ and $Y$ be compact metric spaces and $f : X \to Y$ a
continuous surjection. Assume the $Y$ has trivial reduced K-theory
and that $\Ext\left(K^1(Y), \mathbb Z\right) = 0$. Then the map
$f_* : K_0(X) \to K_0(Y)$ is surjective.
\begin{proof} Let
 $$
X = X^n_1 \sqcup X^n_2 \sqcup \dots \sqcup X^n_{m_n},
 $$
$n \in \mathbb N$, be a  sequence of partitions of $X$ into
non-empty closed and open subsets such that
\begin{enumerate}
\item[i)] the $n+1$'st partition is a refinement of the $n$'th partition,
\item[ii)] $m_{n+1} \leq m_n +1$,
\item[iii)] $C(X,\mathbb Z) = \bigcup_{n=1}^{\infty} \mathcal A_n$,
where $\mathcal A_n$ is the subgroup consisting of the continuous
$\mathbb Z$-valued functions on $X$ that are constant on each
$X^n_i$.
\end{enumerate}
Let $\mathcal B_n$ denote the subgroup of $C(Y,\mathbb Z)$
consisting of the continuous $\mathbb Z$-valued functions on $X$
that are constant on each $f\left(X^n_i\right)$. Then $\mathcal
B_n \subseteq \mathcal B_{n+1}$ and $C(Y,\mathbb Z) =
\bigcup_{n=1}^{\infty} \mathcal B_n$. Note that $\mathcal A_n =
\mathbb Z^{m_n}$ and $\mathcal B_n = \mathbb Z^{k_n}$, when $k_n
\leq m_n$ is the number of elements in the partition of $Y$ which
consists of unions of the $f\left(X^n_i\right)$'s. Hence we also
have $\Hom\left(\mathcal A_n, \mathbb Z\right)  = \mathbb Z^{m_n}$
and $\Hom\left(\mathcal B_n, \mathbb Z\right)  = \mathbb Z^{k_n}$.
Since the map $\Hom\left(\mathcal A_{n+1}, \mathbb Z\right) \to
\Hom\left(\mathcal A_{n}, \mathbb Z\right)$ is surjective for each
$n$, there is an identification $\Hom\left( C(X,\mathbb Z),
\mathbb Z\right)= \varprojlim \Hom \left( \mathcal A_n , \mathbb Z
\right)$. Similarly, there is a also an isomorphism $\Hom\left(
C(Y,\mathbb Z), \mathbb Z\right)= \varprojlim \Hom \left( \mathcal
B_n , \mathbb Z \right)$. Let $\gamma_X : K_0(X) \to \Hom\left(
K_0\left(C(X)\right), \mathbb Z\right)$ and $\gamma_Y : K_0(Y) \to
\Hom\left( K_0\left(C(Y)\right), \mathbb Z\right)$ be the
canonical maps. Combined with the above inverse limit
decompositions we get a commutative diagram
\begin{equation}
\begin{xymatrix}{
{K_0\left(X\right)}  \ar[d]_-{f_*} \ar[r]^-{\gamma_X} & {\Hom \left( K_0\left(C(X)\right), \mathbb Z\right)}  \ar[d]_-{f_*} \ar[r] & {\varprojlim \Hom \left( \mathcal A_n , \mathbb Z \right) }  \ar[d]_-{f_*}  \\
{K_0\left(Y\right)}   \ar[r]^-{\gamma_Y} & \Hom \left( K_0\left(C(Y)\right), \mathbb Z\right)  \ar[r] & {\varprojlim \Hom \left( \mathcal B_n , \mathbb Z \right)} .  }
\end{xymatrix}
\end{equation}
It follows from the universal coefficient theorem of Rosenberg and
Schochet, \cite{RS}, that $\gamma_Y$ is an isomorphism since
$\Ext\left(K^1(Y), \mathbb Z\right) = 0$ by assumption. It follows
that the composition of the maps on the lower row is an
isomorphism because $Y$ has trivial reduced K-theory by
assumption. In the upper row the first map $\gamma_X$ is
surjective by \cite{RS} and the second is surjective because
$C(X,\mathbb Z)$ is a direct summand of $K_0\left(C(X)\right)$. It
suffices therefore to show that
\begin{equation}\label{ann17}
f_* : \varprojlim \Hom\left(\mathcal A_n, \mathbb Z\right) \to \varprojlim \Hom\left(\mathcal B_n, \mathbb Z\right)
\end{equation}
is surjective. To this end we consider a square
\begin{equation}\label{ann1}
\begin{xymatrix}{
{\Hom \left(\mathcal A_n, \mathbb Z\right)}  \ar[d]_-{f_*} & {\Hom \left( \mathcal A_{n+1}, \mathbb Z\right)} \ar[l]_{{i_n}^*} \ar[d]^-{f_*} \\
{\Hom \left( \mathcal B_n, \mathbb Z\right)} & {\Hom \left(\mathcal B_{n+1}, \mathbb Z\right)}  \ar[l]^{{j_n}^*}, }
\end{xymatrix}
\end{equation}
where $i_n : \mathcal A_{n} \to \mathcal A_{n+1}$ and $j_n :
\mathcal B_{n} \to \mathcal B_{n+1}$ are the inclusions. To prove
surjectivity of (\ref{ann17}) it suffices to show that if we are
given $u \in \Hom \left(\mathcal A_n, \mathbb Z\right)$ and $v \in
\Hom \left( \mathcal B_{n+1}, \mathbb Z\right)$ such that
${j_n}^*(v) = f_*(u)$, there is an element $z \in  {\Hom \left(
\mathcal A_{n+1}, \mathbb Z\right)}$ such that ${i_n}^*(z) = u$
and $f_*(z) = v$. This is trivial when $ m_n = m_{n+1}$ because
${i_n}^*$ and $ {j_n}^*$ are identity maps in this case. So assume
that $m_{n+1} = m_n +1$. Let $c_1,c_2, \dots, c_{m_n}$  be the
elements of the partition $\left\{X^{n}_i\right\}_{i=1}^{m_{n}}$
and $c_0',c_1',c_2,c_3, \dots, c_{m_n}$ the elements of the
partition $\left\{X^{n+1}_i\right\}_{i=1}^{m_{n+1}}$, so that
$c_0'$ and $c_1'$ are the only new elements, obtained by a
splitting of $c_1$. We divide the considerations into the
following cases:

\smallskip

$k_{n+1} = k_n +1$: In this case the diagram (\ref{ann1}) takes the form
\begin{equation}\label{ann2}
\begin{xymatrix}{
{\mathbb Z \oplus \mathbb Z^{m_n-1}}  \ar[d]_-{A} & {\mathbb Z^2 \oplus \mathbb Z^{m_n-1}} \ar[l]_{C} \ar[d]^-{D} \\
{\mathbb Z \oplus \mathbb Z^{k_n-1}}   & {\mathbb Z^2 \oplus \mathbb Z^{k_n-1}}  \ar[l]^{B} }
\end{xymatrix}
\end{equation}
where $A,B,C$ and $D$ are matrices of zeroes and ones such that
every column contains one and only one non-zero entry. In the
present case the matrices take the form
 $$
A = \left( \begin{matrix} 1 & \epsilon_2 & \epsilon_3 & \hdots & \epsilon_{m_n} \\ 0 & x_{22}  & x_{23}  & \hdots & x_{2 m_n} \\
0 & x_{32}  & x_{33} & \hdots & x_{3 m_n} \\
\vdots & \vdots & \vdots & \ddots & \vdots \\
0 & x_{k_n,2} &  x_{k_n,3} & \hdots & x_{k_n,m_n} \end{matrix} \right),
 $$

 $$
B = \left( \begin{matrix} 1 & 1 & 0 & 0 & \hdots & 0 \\ 0 & 0  & 1  & 0 & \hdots & 0 \\
0 & 0  & 0 & 1 & \hdots & 0 \\
\vdots & \vdots & \vdots & \vdots & \ddots & \vdots \\
0 & 0 &  0 & 0& \hdots & 1 \end{matrix} \right),
 $$
 $$
C =  \left( \begin{matrix} 1 & 1 & 0 & 0 & \hdots & 0 \\ 0 & 0  & 1  & 0 & \hdots & 0 \\
0 & 0  & 0 & 1 & \hdots & 0 \\
\vdots & \vdots & \vdots & \vdots & \ddots & \vdots \\
0 & 0 &  0 & 0& \hdots & 1 \end{matrix} \right),
 $$
and
 $$
D = \left( \begin{matrix} 1 & 0 & a_2 & a_3 & \hdots & a_{m_n} \\ 0 & 1 & b_2 & b_3 & \hdots & b_{m_n} \\
0 & 0 & x_{22} & x_{23} & \hdots & x_{2 m_n} \\ 0 & 0 & x_{32} &  x_{33} & \hdots & x_{3 m_n} \\
\vdots & \vdots & \vdots & \vdots & \ddots & \vdots \\
0 & 0 & x_{k_n, 2} & x_{k_n,3} & \hdots & x_{k_n, m_n} \end{matrix} \right),
 $$
where $a_i = b_i = 0$, unless $\epsilon_i =1$ in which case
$(a_i,b_i) = (1,0)$ or $(a_i,b_i) = (0,1)$. To find the desired $z
\in \mathbb Z^{m_n+1}$, write $u = \left(u_1,u_2, \dots,
u_{m_n}\right)$ and $v = \left( v_1,v_2, \dots, v_{k_n+1}\right)$.
Then
 $$
z = \left( v_1 - \sum_{i=2}^{m_n} a_iu_i, v_2 - \sum_{i=2}^{m_n}b_i u_i, u_2,u_3, \dots, u_{m_n}\right)
 $$
has the right properties. $z$ is unique in this case.

\smallskip

\smallskip

$k_{n+1} = k_n$: In this case diagram (\ref{ann1}) takes the form
\begin{equation}\label{an32}
\begin{xymatrix}{
{\mathbb Z \oplus \mathbb Z^{m_n-1}}  \ar[d]_-{A} & {\mathbb Z^2 \oplus \mathbb Z^{m_n-1}} \ar[l]_{C} \ar[d]^-{D} \\
{\mathbb Z \oplus \mathbb Z^{k_n-1}}   & {\mathbb Z \oplus \mathbb Z^{k_n-1}}  \ar[l]^{B} }
\end{xymatrix}
\end{equation}
where $A$ and $C$ are as before, but $B$ and $D$ have changed to the identity matrix and
 $$
D = \left( \begin{matrix} 1 & 1 & \epsilon_2 & \epsilon_3 & \hdots & \epsilon_{m_n} \\
0 & 0 & x_{22} & x_{23} & \hdots & x_{2 m_n} \\ 0 & 0 & x_{32} &  x_{33} & \hdots & x_{3 m_n} \\
\vdots & \vdots & \vdots & \vdots & \ddots & \vdots \\
0 & 0 & x_{k_n, 2} & x_{k_n,3} & \hdots & x_{k_n, m_n} \end{matrix} \right),
 $$
respectively. In this case the solution is not unique; if
$\alpha,\beta \in \mathbb Z$ satisfy that $\alpha + \beta = u_1$,
we can use $z = \left(\alpha,\beta, u_2,u_3, \dots,
u_{m_n}\right)$.

\end{proof}
\end{thm}

As a first step in the proof of Theorem \ref{introthm} from the
Introduction, we prove the following:

\begin{thm}\label{THM1}
Let $M,N_1,N_2, N_3, \dots, N_k$ be bounded normal operators such
that $N_iN_j = N_jN_i$ for all $i,j$, and let $F$ be a continuous
function from the joint spectrum of the $N_i$'s onto the spectrum
of $M$ such that
\begin{equation}\label{ass}
F\left(N_1,N_2, \dots, N_k\right) - M \in \mathbb K .
\end{equation}
Assume that the spectrum of $M$ contains no isolated eigenvalue of
finite multiplicity. There are then normcontinuous paths $N^t_i, i
= 1,2, \dots, k, \ t \in [1,\infty)$, of bounded normal operators
such that $N_i^{t}{N_j^{t}} = {N_j^{t}}N_i^{t}$, $N_i-N_i^{t} \in
\mathbb K$ for all $i,j, t$, and
 $$
\lim_{t \to \infty} \left\|F\left(N_1^{t},N_2^{t}, \dots, N_k^{t}\right) - M \right\| = 0.
 $$
\begin{proof} Let $X_0$ be the joint spectrum of $N_1,N_2, \dots, N_k$, and $X$ the joint spectrum of $q(N_1),q(N_2), \dots, q(N_k)$. Then $X$ is a closed subset of $X_0$ and $X/X_0$ is totally disconnected. It follows therefore from Lemma (6.4) of \cite{BDF2} that the extension of $C(X)$ given by the $q\left(N_i\right)$'s is split. Consequently there are commuting normal operators $N_1',N_2',\dots, N_k'$ such that $N_i - N_i' \in \mathbb K$ for all $i$, and the joint spectrum of $N_1',N_2', \dots, N_k'$ is $X$. We may therefore assume from the beginning that the joint spectrum of $N_1,N_2, \dots , N_k$ is equal to the joint spectrum of $q(N_1),q(N_2), \dots, q(N_k)$. Let $\varphi_0 : C\left(\sigma (M)\right) \to \mathbb L(H)$ and $\varphi : C(X) \to \mathbb L(H)$ be the unital $*$-homomorphisms coming from the spectral theory of $M$ and $N_1,N_2, \dots, N_k$, respectively. With $\sigma(M)$ in the role of $Y$ and $F$ in the role of $f$, we are in the setting of Section \ref{K-hom}. By combining the six-term exact sequence of Theorem \ref{cor1} with Theorem \ref{ann6} above we conclude that the pair $(\varphi, \varphi_0)$ represents zero in $\Ext_{\sigma(M),F}(X)$. The desired path of normal operators arise then from condition ii) of Theorem \ref{triv} in the obvious way.
\end{proof}
\end{thm}

\begin{remark}\label{DEremark}
As observed in the proof of Theorem \ref{THM1}, it is
straightforward to reduce the theorem, by use of \cite{BDF2}, to
the case where the joint spectrum of the $N_i$'s is equal to the
joint spectrum $X$ of the $q\left(N_i\right)$'s. After this
reduction the assumption (\ref{ass}) is exactly that the
$*$-homomorphisms $\varphi_+ : C(\sigma(M)) \to \mathbb L(H)$ and
$\varphi_- : C(\sigma(M)) \to \mathbb L(H)$ arising by spectral
theory from $F\left(N_1,N_2, \dots, N_k\right)$ and $M$,
respectively, is a Cuntz-pair, i.e. define an element of
$KK\left(C(\sigma(M)),\mathbb K\right)$. In the case where this
element is trivial in $KK\left(C(\sigma(M)),\mathbb K\right)$ the
conclusion of the theorem follows from Theorem 3.12 of \cite{DE2}.
On the other hand, it is clear that the conclusion of Theorem
\ref{THM1} implies that the element $[\varphi_+,\varphi_-] \in
KK\left(C(\sigma(M)),\mathbb K\right)$ is in the range of $F_* :
KK\left(C(X),\mathbb K\right) \to KK\left(C(\sigma(M)),\mathbb
K\right)$, cf. Lemma \ref{iso6}. What the proof of Theorem
\ref{THM1} does, is to show that this condition is also
sufficient, and always satisfied.
\end{remark}

We can now give the proof of Theorem \ref{introthm}:

\begin{proof}
By spectral theory there is a finite rank projection $P$ such that
$PM = MP$, and a normal operator $M' \in \mathbb
L\left(P^{\perp}H\right)$ such that $\left\|M' -
M|_{P^{\perp}H}\right\| \leq \epsilon, M|_{P^{\perp}H} -M' \in
\mathbb K\left(P^{\perp}H\right)$ and $\sigma(M') =
\sigma_{ess}(M')$. Set $N_i' = P^{\perp}N_i|_{P^{\perp}H} \in
\mathbb L\left(P^{\perp}H\right)$. As in the proof of Theorem
\ref{THM1} we let $X$ denote the joint spectrum of
$q\left(N_1\right),\left(N_2\right), \dots , q\left(N_k\right)$.
Both the $N_i$'s and the $N_i'$'s define an extension of $\mathbb
K$ by $C(X)$ in the sense of \cite{BDF2}, and as argued in the
proof of Theorem \ref{THM1} the extension arising from the $N_i$'s
is split. It follows therefore from Theorem (4.3) of \cite{BDF2}
that the same is true of the extension arising from the $N_i'$'s.
This means that there are commuting normal operators $D_i, i =
1,2, \dots, k$, acting on $P^{\perp}H$ such that $N_i' - D_i \in
\mathbb K\left(P^{\perp}H\right)$ for all $i$, and such that the
joint spectrum of the $D_i$'s is $X$. It follows from the
conditions on $F$ that $F(X) = \sigma_{\ess}(M)$. Since
$\sigma(M') = \sigma_{\ess}(M') = \sigma_{\ess}(M)$ and
$F\left(D_1,D_2, \dots , D_k\right) - M' \in \mathbb
K\left(P^{\perp}H\right)$, it follows from Theorem \ref{THM1} that
there are commuting normal operators $D^{\epsilon}_i$ on
$P^{\perp}H$ such that $D_i -D^{\epsilon}_i \in \mathbb
K\left(P^{\perp}H\right)$ for all $i$ and
$\left\|F\left(D^{\epsilon}_1, D^{\epsilon}_2, \dots,
D^{\epsilon}_k\right) - M' \right\| \leq \epsilon$. Let
$\mu_1,\mu_2, \dots , \mu_L$ be the eigenvalues of $MP$ on $PH$,
each repeated according to its multiplicity so that $L$ is the
rank of $P$. Let $e_1,e_2, \dots, e_L$ be the corresponding
one-dimensional eigenprojections. Since $F$ is surjective by
assumption, there is a $L \times k$ complex matrix
$\left(a_{ij}\right)$ such that $F\left(a_{i1},a_{i2}, \dots ,
a_{ik} \right) = \mu_i$ for all $i$, and we set
 $$
N^{\epsilon}_j =   D_j^{\epsilon}P^{\perp} + \sum_{i=1}^L a_{ij}e_i  .
 $$
Then $N_j - N_j^{\epsilon} \in \mathbb K$ and $\left\| F
\left(N^{\epsilon}_1,N^{\epsilon}_2, \dots , N^{\epsilon}_k
\right) - M \right\| \leq 2\epsilon$.
\end{proof}

We want to point out that the approximation aspect in Theorem
\ref{introthm} and Theorem \ref{THM1}, and hence also in the
theorems of Section \ref{K-hom}, is inevitable. Specifically, we
want to show that in general it is not possible, in the setting of
Theorem \ref{THM1} to find commuting normal operators
$N_1^0,N_2^0, \dots, N^0_k$ such that each $N^0_i$ is a compact
perturbation of $N_i$ and $F\left(N_1^0,N_2^0, \dots, N_k^0\right)
=M$; not even when $k =1$. To this end, let $e_i, i \in \mathbb
Z$, be an orthonormal basis in $H$. Let $z_i, i \in \mathbb Z$, be
a dense sequence in $\mathbb T$ such that $\re z_i \neq \re z_j$
when $i \neq j$ and $\lim_{i \to \infty} z_i = \frac{1}{\sqrt{2}}
+ \frac{i}{\sqrt{2}}$. Define $D \in \mathbb L(H)$ such that $De_i
= 2 \re z_i e_i$ for all $i$, and $T \in \mathbb L(H)$ such that
$Te_i = z_ie_i, i \leq 0$, while $Te_i = \frac{1}{\sqrt{2}}e_i +
\frac{i}{\sqrt{2}}e_{i+1}$ and $Te_{i+1} = \frac{i}{\sqrt{2}}e_i +
\frac{1}{\sqrt{2}}e_{i+1}$ when $i \geq 1$ is odd. Then $T$ is
unitary with $\sigma(T) = \sigma_{ess}(T) = \mathbb T$ and $T^* +
T - D \in \mathbb K$. Since any normal operator $N$ with $N^* + N
= D$ must be diagonal with respect to the basis
$\left\{e_i\right\}_{i \in \mathbb Z}$, such an $N$ can not be a
compact perturbation of $T$.


\begin{thebibliography}{99999}

%\bibitem[A]{A} W. Arveson, {\em Notes on extensions of $C^*$-algebras}, Duke Math. J. {\bf 44} (1977), 329-355.


%\bibitem[Bo]{Bo} F. Boca, {\em Free products of completely positive maps and spectral sets}, J. Funct. Anal. {\bf 97} (1991), 251--263.


%\bibitem[Bl]{Bl} B. Blackadar, {\em K-theory for Operator Algebras}, Springer Verlag, New York, 1986.

%\bibitem[Br]{Br} L. Brown, {\em Ext of certain free product $C^*$-algebras}, J. Operator Theory {\bf 6} (1981), 135-141.


\bibitem[BDF1]{BDF1} L. Brown, R.G. Douglas, P.A. Fillmore {\em Unitary equivalence modulo the compact operators and extensions of $C^*$-algebras}, Proc. Conf. on Operator Theory, Springer Lecture Notes in Math. {\bf 345} (1973), 58-128.



\bibitem[BDF2]{BDF2} \bysame, {\em Extensions of $C^*$-algebras and $K$-homology}, Ann. Math. {\bf 105} (1977), 265-324.


%\bibitem[Cu]{Cu} J. Cuntz, {\em A New Look at KK-theory}, K-theory {\bf 1} (1987), 31-51.



%\bibitem[C]{C} J. Cuntz, {\em The K-groups for free products of $C^*$-algebras}, Proc. Symp. Pure Math. {\bf 38}, Part 1, A.M.S. (1982).


\bibitem[CH]{CH} J. Cuntz and N. Higson, {\em Kuiper's Theorem for Hilbert Modules}, Contemporary Mathematics {\bf 62}, 429-435, AMS, Providence, 1987.

\bibitem[CS]{CS} J. Cuntz and G. Skandalis, {\em Mapping cones and exact sequences in KK-theory}, J. Operator Theory {\bf 15} (1986), 163-180.



%\bibitem[DE]{DE} M. Dadarlat and S. Eilers, {\em On the classification of nuclear $C^*$-algebras}, Proc. London Math. Soc. (3) {\bf 85} (2002), 168-210.


\bibitem[DE]{DE2} \bysame, {\em Asymptotic unitary equivalence in KK-theory}, K-theory {\bf 23} (2001), 305-322.


%\bibitem[G1]{G1} E. Germain, {\em $KK$-theory of reduced free-product $C\sp *$-algebras}, Duke Math. J. {\bf 82} (1996), 707-723.


%\bibitem[G2]{G2} \bysame , {\em KK-theory of the full free product of unital $C^*$-algebras}, J. Reine angew. Math. {\bf 485} (1997),1-10.


%\bibitem[H]{H} N. Higson, {\em $C^*$-Algebra Extension Theory and Duality}, J. Funct. Anal. {\bf 129} (1995),349-363.



\bibitem[K1]{K1} G.G. Kasparov, {\em  The operator K-functor and extensions of $C^*$-algebras}, Izv. Akad. Nauk. SSSR, Ser. Math. {\bf 44} (1980), 571-636; Math. USSR Izvestija {\bf 16}, 513-572.



\bibitem[K2]{K2} \bysame, {\em Hilbert $C^*$-modules: theorems of Stinespring and Voiculescu}, J. Operator Theory {\bf 4} (1980), 133-150.




\bibitem[MT]{MT} V. Manuilov, K. Thomsen, {\em The group of unital $C^*$-extensions},  Preprint, 2005.

%\bibitem[MT]{MT} V. Manuilov, K. Thomsen, {\em $E$-theory is a special case of $KK$-theory},  Proc. London Math. Soc. (3)  88  (2004),  no. 2, 455--478.

\bibitem[M]{M} J.A. Mingo, {\em On the contractibility of the unitary group of the Hilbert space over a $C^*$-algebra}, Integral Eq. Operator Theory {\bf 5} (1982), 888-891.

%\bibitem[OP]{OP} C. Olsen and G.K. Pedersen, {\em Corona $C^*$-algebras and their applications to lifting problems}, Math. Scand. {\bf 64} (1989), 63-86.



\bibitem[Pa]{Pa} W.L. Paschke, {\em K-theory for commutants in the Calkin algebra}, Pacific J. Math. {\bf 95} (1981), 427-434.



\bibitem[RS]{RS} J. Rosenberg and C. Schochet, {\em The K\"unneth theorem and the universal coefficient theorem for Kasparov's generalized K-functor}, Duke Math. {\bf 55} (1987), 431-474.




\bibitem[S]{S} G. Skandalis, {\em Kasparov's bivariant $K$-theory and applications},  Exposition. Math.  9  (1991),  no. 3, 193-250.



%\bibitem[Th0]{Th0}  K. Thomsen,{\it Homotopy classes of $*$-homomorphisms between stable $C\sp *$-algebras and their multiplier algebras},  Duke Math. J. {\bf 61} (1990), 67-104.

\bibitem[Th1]{Th1}  K. Thomsen, {\em On absorbing extensions}, Proc. Amer. Math. Soc. {\bf 129} (2001), 1409-1417.

\bibitem[Th2]{Th2} \bysame, {\em Homotopy invariance in $E$-theory}, Preprint.

\bibitem[Th3]{Th3} \bysame, {\em On the KK-theory and the E-theory of amalgamated free products of $C^*$-algebras}, J. Func. Anal. {\bf 201} (2003), 30-56.


%\bibitem[Th4]{Th4} \bysame, {\em Limits of certain subhomogeneous $C^*$-algebras} M\'emoires de la Soci\'et\'e Math\'e matique de France {\bf 71} (1997).



%\bibitem[V]{V} D. Voiculescu,{\em A non-commutative Weyl-von Neumann theorem}, Rev. Roum. Math. Pures et Appl. {\bf 21} (1976), 97-113.

\end{thebibliography}
\end{document}